# QUEUEING SYSTEMS WITH MANY SERVERS: NULL CONTROLLABILITY IN HEAVY TRAFFIC[1]

By Rami Atar, Avi Mandelbaum and Gennady Shaikhet

*Technion—Israel Institute of Technology*

A queueing model has $J \geq 2$ heterogeneous service stations, each consisting of many independent servers with identical capabilities. Customers of $I \geq 2$ classes can be served at these stations at different rates, that depend on both the class and the station. A system administrator dynamically controls scheduling and routing. We study this model in the central limit theorem (or heavy traffic) regime proposed by Halfin and Whitt. We derive a diffusion model on $\mathbb{R}^I$ with a singular control term that describes the scaling limit of the queueing model. The singular term may be used to constrain the diffusion to lie in certain subsets of $\mathbb{R}^I$ at all times $t > 0$. We say that the diffusion is *null-controllable* if it can be constrained to $\mathbb{X}_-$, the minimal closed subset of $\mathbb{R}^I$ containing all states of the prelimit queueing model for which all queues are empty. We give sufficient conditions for null controllability of the diffusion. Under these conditions we also show that an analogous, asymptotic result holds for the queueing model, by constructing control policies under which, for any given $0 < \varepsilon < T < \infty$, all queues in the system are kept empty on the time interval $[\varepsilon, T]$, with probability approaching one. This introduces a new, unusual heavy traffic "behavior": On one hand, the system is critically loaded, in the sense that an increase in any of the external arrival rates at the "fluid level" results with an overloaded system. On the other hand, as far as queue lengths are concerned, the system behaves as if it is underloaded.

**1. Introduction.** We consider a multiclass queueing model with heterogeneous service stations, each consisting of many independent servers with identical capabilities. The servers offer service to different classes of customers at rates that may depend on both the station and the class. A sys-

Received July 2005; revised March 2006.

[1]Supported in part by the Israel Science Foundation Grant 126/02 and the Technion funds for the promotion of research and sponsored research.

tem administrator dynamically controls all scheduling and routing in the system. The model is considered in the heavy traffic parametric regime, first proposed by Halfin and Whitt [7], in which the number of servers at each station and the arrival rates grow without bound, while keeping, in an appropriate sense, a critically balanced system. Both the model and the parametric regime have recently received much attention, especially in relation to large telephone call centers (see [6] and references therein).

When studying queueing models in heavy traffic, one considers a sequence of models parametrized by $n \in \mathbb{N}$ that, under a law of large numbers (LLN) limit, give rise to a *fluid model* which is critically loaded. Typically, an attempt is then made to prove that appropriately scaled fluctuations of the queueing model about the fluid model converge to a diffusion. If, as in the current setting, a control problem is associated to the queueing model, then a similar approach gives rise to a *controlled diffusion model*. In this case, a natural goal is to prove that, given a cost criterion, the (suitably scaled) value of a control problem for the queueing model converges to one for the diffusion model. Moreover, it is often the case that solving the diffusion model for optimal controls helps understand how to construct control schemes for the queueing model that are asymptotically optimal. A similar approach is taken in the current paper. However, rather than a cost criterion, our formulation will be concerned with a certain property observed for the diffusion model and shown to be inherited, in an asymptotic sense, by the queueing model. The property, which is unusual in heavy traffic formulations, will have to do with the ability to maintain empty queues.

We let $I$ denote the number of customer classes, and let $X^n(t)$, $t \geq 0$, denote an $I$-dimensional process for which the $i$th component $X_i^n(t)$ represents the total number of class-$i$ customers in the system at time $t$, in the $n$th system. The fluctuations alluded to above are denoted by $\hat{X}^n(t)$, $t \geq 0$. These fluctuations are scaled in such a way that $\hat{X}^n$ gives rise, as weak limits are taken formally, to a (controlled) diffusion process denoted by $X(t)$, $t \geq 0$. Moreover, one has that

$$\hat{X}^n(t) \in \mathbb{X}_- := \left\{ x \in \mathbb{R}^I : \sum_{i=1}^I x_i \leq 0 \right\}$$

holds if and only if the total number of customers in the $n$th system at time $t$ is less than or equal to the total number of servers. Ideally, if one could freely rearrange customers in the system without any constraints, it follows that one could maintain empty queues whenever $\hat{X}^n(t) \in \mathbb{X}_-$. We will thus refer to $\mathbb{X}_-$ as the *null domain*. Although the queueing model considered in this paper is subject to additional constraints (e.g., a station may offer service to only some of the classes), the null domain will play the same role in the asymptotic regime under consideration.

An important feature of the controlled diffusion model derived in this paper is that the stochastic differential equation describing it has a *singular control* term, that is, a control process with sample paths that are locally of bounded variation, the increments of which take values in a fixed cone $\mathbb{C}$ of $\mathbb{R}^I$. The singular term may be used to constrain the diffusion to lie in certain subsets of $\mathbb{R}^I$ at all times $t > 0$. We say that the diffusion is *null-controllable* if

$$(1) \qquad \mathbb{C} \cap \mathbb{X}_-^o \neq \varnothing,$$

where $\mathbb{X}_-^o$ denotes the interior of $\mathbb{X}_-$, because, under this condition, the diffusion can be constrained to $\mathbb{X}_-$. Condition (1) will be given in explicit form in terms of the model parameters [see (33)]. Our main result shows that under (1) one can construct control policies for the queueing model in such a way that, for every given $0 < \varepsilon < T < \infty$, all queues are kept empty on the time interval $[\varepsilon, T]$ with probability approaching 1, as $n \to \infty$. We will refer to such behavior as *asymptotic null-controllability*. We will, in fact, consider two versions of the problem: One, referred to as *preemptive scheduling*, in which service to a customer can be interrupted and resumed at a later time (possibly in a different station). The other, referred to as *nonpreemptive scheduling*, where service to customers may not be interrupted before service is completed. The treatment of the nonpreemptive case is more complicated than that of the preemptive case. Thus, to keep the exposition simple, we have limited ourselves in the nonpreemptive case to the simplest possible network structure where null-controllability can show up: two customer classes and two service stations.

Our results on asymptotic null-controllability can be regarded as the demonstration of a new, unusual heavy traffic behavior. On one hand, the system is critically loaded. Indeed, as intuitively expected (and precisely stated in Proposition 2.1), an increase in any of the external arrival rates at the fluid level results with an overloaded system, in the sense that large queues inevitably build up. On the other hand, the system behaves as if it is underloaded as far as the capability of maintaining empty queues is concerned.

Singular control arises in connection with queueing systems in heavy traffic in many references. The singular term is often associated with positivity or finite buffer constraints for the queue length process (see [9], Chapter 8, for discussion and further references), with admission control (ibid.) or with constraints on the so-called workload process to lie in a given cone [5]. The source for the singular term in the current setting is, however, quite different, and it has to do with the fact that a many-server limit is taken. To explain this point, consider a network in which customers of classes, say, 1 and 2, can be served at both station $A$ and station $B$ (these two classes and two stations could be a subset of a larger setup). Assume that the network operates under preemptive scheduling. Suppose that at a certain moment one selects, say, $r$ class-1 customers that are in service at station $A$ and $r$ class-2 customers in station $B$ and considers the option of interchanging their position, so that the $r$ class-1 (class-2) customers that were selected are moved to station $B$ (resp., $A$). Since the service rates may depend on both the class and the station, the average rate at which components 1 and 2 of $X^n$ change at that moment may be different depending on whether such an interchange takes place or not. If the interchange is performed instantaneously, the rates alluded to above will change abruptly. Moreover, since in both stations the number of servers is assumed to be large, this effect can be amplified by letting $r$ be large. This, as scaling limit is taken, results with a control term that can have arbitrarily large increments over a given time interval. An appropriate formulation for such a controlled diffusion model will thus have a singular control term. Next, let $\mathcal{G}$ denote the graph with classes and stations as vertices, and with $(i,j)$ pairs as edges if and only if station $j$ can serve class $i$. It is instructive to note that our explanation above relies on a certain property of the graph $\mathcal{G}$, namely, that it contains a cycle with vertices 1, 2, $A$ and $B$. Indeed, this is just another way of saying that customers of each of the two classes can be served in both stations, as we have assumed. A heuristic similar to the one described above will lead to a singular term in the diffusion model whenever $\mathcal{G}$ contains cycles, whether with four vertices or more. Thus, in general, cycles contained in $\mathcal{G}$ will play an important role in the singular control formulation of the diffusion model.

Note that the phenomenon described above is indeed a result of the many server setting, because it must be possible to ocasionally let the number $r$ referred to above take large values, and $r$ is clearly limited by the number of servers. In particular, this phenomenon is not seen in what is sometimes referred to as "conventional" heavy traffic, where diffusion limits are obtained for systems with a fixed number of servers (and large service rates; cf. [11]).

Recall that the so-called fluid model describes the LLN limit of basic quantities of the queueing model. One ingredient of the fluid model is a (deterministic, constant) matrix denoted by $\xi^*$. The entry $\xi^*_{ij}$ represents the (large $n$ limit) fraction of the number of servers in station $j$ that serve class-$i$ customers. One refers to $(i,j)$ pairs that are edges of $\mathcal{G}$ as *activities*, and to activities $(i,j)$ for which $\xi^*_{ij} > 0$ as *basic* activities. In particular, the number of servers that are engaged in nonbasic activities ($\xi^*_{ij} = 0$) is of order $o(n)$ (where the total number of servers is proportional to $n$), in a sense that can be made precise (see, e.g., Lemma 4.1). The policy that we shall propose for the preemptive case will basically mimic the construction of a constrained diffusion. Namely, a special rearrangement of customers in the service stations will take place whenever the process $\hat{X}^n$ reaches close to the boundary of $\mathbb{X}_-$ (from the inside). In this rearrangement, the number of servers allocated to work in a certain nonbasic activity will be much larger than the typical fluctuations of the process $X^n$ [but still $o(n)$]. This rearrangement will have an effect on the dynamics of $\hat{X}^n$ that is a reminiscent of that of a Skorohod mechanism [10] on a constrained diffusion. Namely, it will constrain $\hat{X}^n$ to $\mathbb{X}_-$ by making the term that, in the limit, shows up as singular control, large. Because of (1), this can be performed in such a way that the term points into the domain. Although this question is not directly addressed in this paper, it is expected, in fact, that the restriction to the time interval $[\varepsilon, T]$ of the processes $\hat{X}^n$ converge to a constrained diffusion as $n \to \infty$. The picture is quite different in the nonpreemptive setting. The arrangement of customers can only be controlled indirectly (via routing decisions), and the constraining mechanism of the diffusion cannot be faithfully mimicked. A convergence result as above is not to be expected; in fact, the processes $\hat{X}^n$ that we construct in this case are not even tight. The main idea behind the policy proposed in this case is to assure that a relatively large portion of the servers are engaged in the nonbasic activity at all times, rather than at times when the boundary is reached as in the preemptive case.

We reiterate that it is the presence of cycles in $\mathcal{G}$ that induces a singular control term in the diffusion. A model similar to the current one is studied in [1] and [2] under structural assumptions that are complementary to those of this paper, in the sense that the graph $\mathcal{G}$ is assumed there to be a tree. Indeed, in these references the diffusion has no singular control term, and a phenomenon as described in this paper does not occur. Finally, we would like to mention that one can also consider a setting in which the diffusion has a singular term as in the current paper, but null controllability does not hold, and approach the model from a control theoretic viewpoint so as to minimize costs of interest. This will be subject for future study.

The organization of the paper is as follows. Some notation is introduced at the end on this section. In Section 2 we first introduce the model and describe how its parameters are rescaled. We then present the main step toward the derivation of the diffusion model in Theorem 2.1, where a representation of the prelimit queueing model is provided. The diffusion model, obtained from this representation under the scaling limit, is stated in equations (30)–(32). We then state that, under (1), the diffusion can be constrained to $\mathbb{X}_-$ (Theorem 2.2), and provide the main results regarding asymptotic null-controllability in the preemptive case (Theorem 2.3) and in the nonpreemptive case (Theorem 2.4). At the end of Section 2 we give numerical examples, and demonstrate that asymptotic null-controllability cannot be expected in overloaded models (Proposition 2.1). Section 3 contains the proofs of Theorems 2.2 and 2.3. The proof of Theorem 2.4 appears in Section 4. The Appendix contains the proofs of Theorem 2.1, Proposition 2.1 and some auxiliary results.

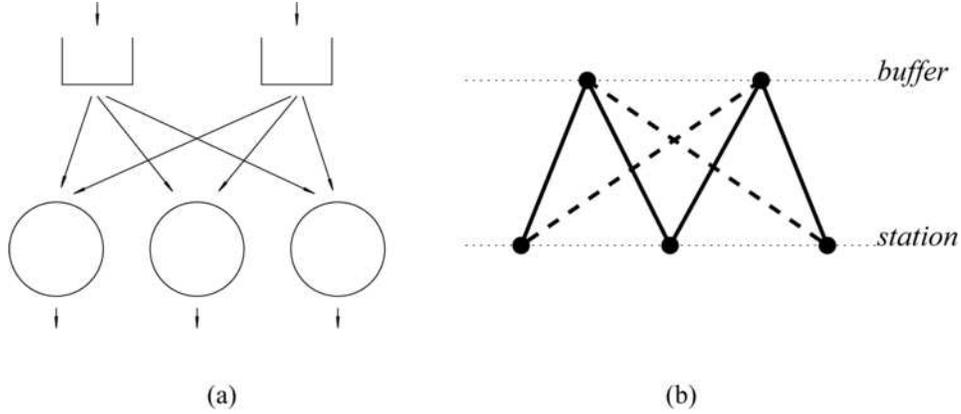

FIG. 1. (a) *A network with* 2 *classes and* 3 *stations.* (b) *Corresponding graph* $\mathcal{G}$ *with basic and nonbasic activities (solid and, resp., dashed lines). The subgraph* $\mathcal{G}_b$ *of* $\mathcal{G}$ *is a tree.*

NOTATION. For $x \in \mathbb{R}^I$, let $\|x\| = \sum_{i=1}^{I} |x_i|$. For two column vectors $v, u$ of the same dimension $v \cdot u$ denotes their scalar product. The symbols $e_i$ denote the unit coordinate vectors and $e = (1, \ldots, 1)'$. The dimension of $e$ may change from one expression to another. For $E$ a metric space, we denote by $\mathbb{D}(E)$ the space of all cadlag functions (i.e., right continuous and having left limits) from $\mathbb{R}_+$ to $E$. We endow $\mathbb{D}(E)$ with the usual Skorohod topology (cf. [4]). If $X^n$, $n \in \mathbb{N}$ and $X$ are processes with sample paths in $\mathbb{D}(E)$, we write $X^n \Rightarrow X$ to denote weak convergence of the measures induced by $X^n$ [on $\mathbb{D}(E)$] to the measure induced by $X$. For a collection $\mathcal{A}$ of random variables, $\sigma\{\mathcal{A}\}$ denotes the sigma-field generated by this collection. If $X$ is an $\mathbb{R}^I$-valued process (or function from $\mathbb{R}_+$ to $\mathbb{R}^I$), $\|X\|_t^* = \sup_{0 \le u \le t} \|X(u)\|$, and if $X$ takes real values, $|X|_t^* = \sup_{0 \le u \le t} |X(u)|$.

## 2. Setting and main results.

2.1. *Queueing model.* A precise description of the queueing model is as follows. A complete probability space $(\Omega, \mathcal{F}, P)$ is given, supporting all stochastic processes defined below. The processes will all be constructed in such a way that they have cadlag sample paths with probability 1. Expectation with respect to $P$ is denoted by $E$. The queueing model is parametrized by $n \in \mathbb{N}$. It has $I \ge 2$ customer classes and $J \ge 2$ service stations [see Figure 1(a)]. Station $j$ has $N_j^n$ identical servers working independently. The classes are labeled as $1, \ldots, I$ and the stations as $I+1, \ldots, I+J$:

$$\mathcal{I} = \{1, \ldots, I\}, \qquad \mathcal{J} = \{I+1, \ldots, I+J\}.$$

Arrivals are modeled as renewal processes with finite second moment for the interarrival time. More precisely, we are given parameters $\lambda_i^n > 0$, $i \in \mathcal{I}$, $n \in \mathbb{N}$, and independent sequences of strictly positive i.i.d. random variables $\{\check{U}_i(k), k \in \mathbb{N}\}$, $i \in \mathcal{I}$, with mean $E\check{U}_i(1) = 1$ and variance $C_{U,i}^2 = \operatorname{Var}(\check{U}_i(1)) \in [0, \infty)$. With $\sum_1^0 = 0$, the number of class-$i$ arrivals up to time $t$ for the $n$th system is given by

$$A_i^n(t) = \sup\left\{l \geq 0 : \sum_{k=1}^{l} \frac{\check{U}_i(k)}{\lambda_i^n} \leq t\right\}, \qquad t \geq 0.$$

For $i \in \mathcal{I}, j \in \mathcal{J}$ and $n \in \mathbb{N}$, we are given parameters $\mu_{ij}^n \geq 0$, representing the service rate of a class-$i$ customer by a server of station $j$. If $\mu_{ij}^n = 0$, we say that class-$i$ customers cannot be served at station $j$. Consider the graph $\mathcal{G}$ having a vertex set $\mathcal{I} \cup \mathcal{J}$ and an edge set

$$\mathcal{E} = \{(i,j) \in \mathcal{I} \times \mathcal{J} : \mu_{ij}^n > 0\}.$$

We assume that $\mathcal{E}$ does not depend on $n$. We denote $i \sim j$ if $(i,j) \in \mathcal{E}$. A class-station pair $(i,j)$ is said to be an *activity* if $i \sim j$, or, equivalently, if class-$i$ customers can be served at station $j$. For every $(i,j) \in \mathcal{I} \times \mathcal{J}$, we denote by $\Psi_{ij}^n(t)$ the number of class-$i$ customers being served in station $j$ at time $t$. By definition,

(2) $$\Psi_{ij}^n(t) = 0 \qquad \text{for } i \not\sim j.$$

Service times are modeled as independent exponential random variables. To this end, let $S_{ij}^n$, $(i,j) \in \mathcal{I} \times \mathcal{J}$ be Poisson processes with rate $\mu_{ij}^n$ (where a Poisson process of zero rate is the zero process), mutually independent and independent of the arrival processes. Note that the time up to $t$ devoted to a class-$i$ customer by a server, summed over all servers of station $j$, is given as $\int_0^t \Psi_{ij}^n(s)\,ds$. The number of service completions of class-$i$ customers by all servers of station $j$ up to time $t$ is, by assumption, given by $S_{ij}^n(\int_0^t \Psi_{ij}^n(s)\,ds)$. The precise description of the processes $\Psi^n = (\Psi_{ij}^n, i \in \mathcal{I}, j \in \mathcal{J})$ is not given at this point. We do emphasize, however, that they will be constructed in such a way that future service completion times are independent of the current state, which results with independent exponential service times (cf. [2]). We note in passing that whereas renewal arrivals are quite natural in the Halfin–Whitt setting, assumptions on service times that go beyond exponential (not to be dealt with in this paper) lead to far more complicated diffusion models [12].

The processes $A^n$ and $S^n$ will be referred to as the *primitive processes*.

Denoting by $X_i^n(t)$ the number of class-$i$ customers in the system (meaning in the queue or being served) at time $t$, and setting $X_i^{0,n} = X_i^n(0)$, it is clear from the above that

(3) $$X_i^n(t) = X_i^{0,n} + A_i^n(t) - \sum_{j \in \mathcal{J}} S_{ij}^n\left(\int_0^t \Psi_{ij}^n(s)\,ds\right), \qquad i \in \mathcal{I}, t \geq 0.$$

For simplicity, the initial conditions $X_i^{0,n}$ are assumed to be deterministic. Finally, we introduce the processes $Y_i^n(t)$, representing the number of class-$i$ customers in the queue (not being served) at time $t$, and $Z_j^n(t)$, representing the number of servers at station $j$ that are idle at time $t$. Clearly, we have the following relations:

$$Y_i^n(t) + \sum_{j \in \mathcal{J}} \Psi_{ij}^n(t) = X_i^n(t), \qquad i \in \mathcal{I}, \tag{4}$$

$$Z_j^n(t) + \sum_{i \in \mathcal{I}} \Psi_{ij}^n(t) = N_j^n, \qquad j \in \mathcal{J}. \tag{5}$$

Also, the following holds by definition:

$$Y_i^n(t) \geq 0, \qquad Z_j^n(t) \geq 0, \qquad \Psi_{ij}^n(t) \geq 0, \qquad i \in \mathcal{I}, j \in \mathcal{J}, t \geq 0. \tag{6}$$

We will write $X^n$ for the vector $(X_i^n, i \in \mathcal{I})$ and similarly $Y^n = (Y_i^n, i \in \mathcal{I})$, $Z^n = (Z_j^n, j \in \mathcal{J})$.

2.2. *Control and rescaling.* Equations (2)–(6) indicate some properties of the processes involved, but they do not characterize these processes, because the control processes $\Psi^n$ have not yet been described. This is the purpose of the following definitions.

*Preemptive scheduling.* We will regard scheduling as preemptive if service to a customer can be stopped and resumed at a later time, possibly in a different station. Formally, such a scheduling is a scheme according to which one selects $\Psi^n(t)$ at every $t$. In this paper we will be concerned only with scheduling of feedback form, in the sense that the selection of $\Psi^n(t)$ depends only on $X^n(t)$, for every $t$. The precise definition is as follows.

DEFINITION 2.1. Let $n$ be given. We say that a map $f^n : \mathbb{Z}_+^I \to \mathbb{Z}_+^{IJ}$ is a preemptive resume scheduling control policy (P-SCP) and $X^n$ is the controlled process corresponding to $f^n$, initial data $X^{0,n}$ and primitive processes $A^n$ and $S^n$, if $\Psi^n(t) = f^n(X^n(t))$ and equations (2)–(6) hold.

*Nonpreemptive scheduling.* By this, we mean that it is impossible to interrupt service to customers. Thus, the quantities $\Psi_{ij}^n$ cannot be directly controlled, but they are affected by the routing decisions according to the following equation:

$$\Psi_{ij}^n(t) = \Psi_{ij}^n(0) + B_{ij}^n(t) - S_{ij}^n\left(\int_0^t \Psi_{ij}^n(s)\, ds\right). \tag{7}$$

Above, for each $(i,j) \in \mathcal{I} \times \mathcal{J}$, $B_{ij}^n$ is a nondecreasing $\mathbb{Z}_+$-valued process starting from zero, that increases by $k$ each time $k$ class-$i$ customers are routed to station $j$. Of course, $B_{ij}^n = 0$ for $i \not\sim j$. In nonpreemptive scheduling, a control policy is a scheme for selecting $B_{ij}^n(t)$ for every $t$. In this paper we will need a randomized formulation, in which the scheme according to which $B_{ij}^n$ are determined may depend on an auxiliary stochastic process. In addition, in our formulation we will only need the routing decisions to take place at the times when arrivals occur. To this end, for $i \in \mathcal{I}$ and $n \in \mathbb{N}$, we let $\tau_k^{n,i}$ denote the time of the $k$th jump of the process $A_i^n$, that is, the time of the $k$th class-$i$ arrival. Finally, it will be assumed that all customers in service at time zero begin their service at that time, and the initial arrangement of these customers in the stations, that is, $\Psi_{ij}^n(0)$, is determined by the policy. We write $\Xi^n(t) = (X^n(t), Y^n(t), Z^n(t), \Psi^n(t))$.

DEFINITION 2.2. Let $n$ be given. We say that a triplet $(\Psi^{0,n}, F^n, R^n)$ is a non-preemptive scheduling control policy (N-SCP) and $X^n$ is the controlled process corresponding to $(\Psi^{0,n}, F^n, R^n)$, initial data $X^{0,n}$ and primitive processes $A^n$, $S^n$ if the following hold:

(i) $\Psi^{0,n} \in \mathbb{Z}_+^{IJ}$ is the initial condition for $\Psi^n$, that is, $\Psi^n(0) = \Psi^{0,n}$. In particular, it satisfies $\Psi_{ij}^{0,n} = 0$ for $i \not\sim j$, $\sum_j \Psi_{ij}^{0,n} \leq X_i^{0,n}$ and $\sum_i \Psi_{ij}^{0,n} \leq N_j^n$.

(ii) $R^n$ is a collection $\{R^{n,i}, i \in \mathcal{I}\}$, where, for each $i$, $R^{n,i} = \{R_k^{n,i} : k \in \mathbb{N}\}$ is a sequence of $\mathbb{R}$-valued independent random variables. The sequences $R^{n,i}$ are mutually independent and independent of the primitive processes.

(iii) $F^n$ is a collection $\{F_{ij}^n, (i,j) \in \mathcal{I} \times \mathcal{J}\}$ of measurable maps $F_{ij}^n : \mathbb{R}_+ \times \mathbb{Z}_+^I \times \mathbb{Z}_+^I \times \mathbb{Z}_+^J \times \mathbb{Z}_+^{IJ} \times \mathbb{R} \to \mathbb{Z}_+$ and for each $(i,j) \in \mathcal{I} \times \mathcal{J}$, $B_{ij}^n$ is given in the form
$$B_{ij}^n(t) = \sum_{k \,:\, \tau_k^{n,i} \leq t} F_{ij}^n(\tau_k^{n,i}, \Xi^n(\tau_k^{n,i}-), R_k^{n,i}), \qquad t \geq 0.$$

(iv) Equations (2)–(7) hold.

REMARK 1. The restriction to policies in which decisions take place only when arrivals occur may not appear to be natural, and one could consider extensions of this definition, say, by allowing decisions to take place upon arrivals or service completions, or even continuously in time. Note, however, that this restriction does not affect our result (i.e., Theorem 2.4) that is concerned with the existence (and construction) of N-SCP with a certain property: Clearly, the existence of such policies under the preceding definition implies the existence of policies under any extention of it.

REMARK 2. Existence and uniqueness of the processes $X^n$ and $\Psi^n$ (given the primitive processes) is easily obtained by induction on the jump times of the primitive processes. In addition, one can argue that the future service completion times are independent of all that has occured up to the current time. For a precise statement and an argument to this effect in a more restriced setup, see Proposition 1 of [3]. This can be adapted to the current setup, but we have omitted the details.

*Fluid scaling.* We assume that the parameters of the processes involved are scaled in the following way. There are constants $\lambda_i, \nu_j \in (0, \infty)$, $i \in \mathcal{I}$, and $\mu_{ij} \in (0, \infty)$, $(i,j) \in \mathcal{E}$ such that

$$(8) \qquad n^{-1}\lambda_i^n \to \lambda_i, \qquad \mu_{ij}^n \to \mu_{ij}, \qquad n^{-1}N_j^n \to \nu_j.$$

We also define $\mu_{ij} = 0$ for $i \not\sim j$. Note that this is consistent with (8) because $\mu_{ij}^n = 0$ for $i \not\sim j$. Let

$$(9) \qquad \bar{\mu}_{ij} = \nu_j \mu_{ij}, \qquad (i,j) \in \mathcal{I} \times \mathcal{J},$$

and consider the following linear program:

$$(10) \quad \begin{aligned} &\text{Minimize } \rho \in \mathbb{R}_+ \text{ subject to} \\ &\sum_{j \in \mathcal{J}} \bar{\mu}_{ij}\xi_{ij} = \lambda_i, \qquad \sum_{i \in \mathcal{I}} \xi_{ij} \leq \rho, \qquad \xi_{ij} \geq 0, \qquad i \in \mathcal{I}, \ j \in \mathcal{J}. \end{aligned}$$

Throughout, we assume that the fluid model is critically loaded. More precisely, we will assume that the *heavy traffic condition* [8] holds: There exists a unique optimal solution $(\xi^*, \rho^*)$ to the linear program (10) and, moreover, $\sum_{i \in \mathcal{I}} \xi_{ij}^* = 1$ for all $i \in \mathcal{J}$ (and, consequently, $\rho^* = 1$). We shall keep the notation $\xi_{ij}^*$ throughout the paper. We also let

$$(11) \qquad x_i^* = \sum_j \xi_{ij}^* \nu_j, \qquad \psi_{ij}^* = \xi_{ij}^* \nu_j, \qquad i \in \mathcal{I}, j \in \mathcal{J},$$

and refer to the quantities $\xi^*, \psi^*, x^*$ as the *static fluid model*, or just *fluid model* for short (see in [2] what these quantities intuitively represent).

Following [8], an activity $(i,j) \in \mathcal{E}$ is said to be *basic* if $\xi_{ij}^* > 0$. Define the graph of basic activities $\mathcal{G}_b$ to be the subgraph of $\mathcal{G}$ having $\mathcal{I} \cup \mathcal{J}$ as a vertex set, and the collection $\mathcal{E}_b$ of basic activities as an edge set. We will also denote the set of nonbasic activities as $\mathcal{E}_{nb} = \mathcal{E} \setminus \mathcal{E}_b$ [see Figure 1(b)].

Like in some other papers that study a similar fluid model in heavy traffic (e.g., [8]), we will have one more assumption in this paper about the fluid model, namely, that the *complete resource pooling* condition holds. This condition expresses, in a sense, a strong mode of cooperation between the service stations. More precisely, one of the equivalent formulations of this condition (see [8]) is that all vertices in $\mathcal{J}$ communicate via edges in $\mathcal{G}_b$. It was shown by Williams [13] that this condition is equivalent to the condition that the basic activities form a tree. Thus, we assume throughout:

$$(12) \qquad \text{The graph } \mathcal{G}_b \text{ is a tree.}$$

We now introduce a rescaled version of the processes describing the queueing model:

$$\bar{X}_i^n(t) = n^{-1} X_i^n(t),$$
$$\bar{Y}_i^n(t) = n^{-1} Y_i^n(t),$$
$$\bar{Z}_j^n(t) = n^{-1} Z_j^n(t),$$
$$\bar{\Psi}_{ij}^n(t) = n^{-1} \Psi_{ij}^n(t).$$

Denote $\bar{X}^n = (\bar{X}_i^n, i \in \mathcal{I})$, and use a similar notation for $\bar{Y}^n$, $\bar{Z}^n$ and $\bar{\Psi}^n$. We will sometimes consider $\bar{X}^n$, $\bar{Y}^n$ and $\bar{Z}^n$ as column vector-valued processes. Heuristically, one expects that $(\bar{X}^n, \bar{Y}^n, \bar{Z}^n, \bar{\Psi}^n) \Rightarrow (x^*, 0, 0, \psi^*)$, and this is indeed the case under appropriate conditions [see, e.g., equation (57) and Lemma 4.1 for such statements in the preemptive and, resp., nonpreemptive case]. For this reason, these processes are referred to as the fluid-level rescaled processes.

*Diffusion scaling.* We further assume that there are constants $\hat{\lambda}_i, \hat{\mu}_{ij} \in \mathbb{R}$, $i \in \mathcal{I}$, $j \in \mathcal{J}$, such that

(13) $\quad \hat{\lambda}_i^n := n^{1/2}(n^{-1}\lambda_i^n - \lambda_i) \to \hat{\lambda}_i, \qquad \hat{\mu}_{ij}^n := n^{1/2}(\mu_{ij}^n - \mu_{ij}) \to \hat{\mu}_{ij},$

(14) $\quad \hat{N}_j^n := n^{1/2}(n^{-1}N_j^n - \nu_j) \to 0.$

We introduce a centered, rescaled version of the primitive processes:

(15) $\quad \hat{A}_i^n(t) = n^{-1/2}(A_i^n(t) - \lambda_i^n t), \qquad \hat{S}_{ij}^n(t) = n^{-1/2}(S_{ij}^n(nt) - n\mu_{ij}^n t).$

Similarly, the processes representing the queueing model are centered about the fluid model quantities and rescaled:

(16) $\qquad \hat{X}_i^n(t) = n^{-1/2}(X_i^n(t) - nx_i^*),$

(17) $\qquad \hat{Y}_i^n(t) = n^{-1/2} Y_i^n(t), \qquad \hat{Z}_j^n(t) = n^{-1/2} Z_j^n(t),$

(18) $\qquad \hat{\Psi}_{ij}^n(t) = n^{-1/2}(\Psi_{ij}^n(t) - \psi_{ij}^* n).$

The processes denoted with hats will be referred to as diffusion-level rescaled processes. Similarly to the fluid-level processes, define $\hat{X}^n = (\hat{X}_i^n, i \in \mathcal{I})$, with an analogous definition for $\hat{Y}^n$, $\hat{Z}^n$ and $\hat{\Psi}^n$, and consider $\hat{X}^n$, $\hat{Y}^n$ and $\hat{Z}^n$ as column vector-valued processes.

*Scaling of initial conditions.* It is assumed that there are constants $x_i$, $i \in \mathcal{I}$, such that the initial conditions satisfy

(19) $\qquad\qquad\qquad \hat{X}_i^{0,n} := \hat{X}_i^n(0) \to x_i.$

Throughout, $x = (x_i, i \in \mathcal{I})$.

2.3. *Main results.* Our first result expresses a relation directly between the diffusion-level processes. Although its proof requires only some elementary manipulations of the relations (2)–(6) and (16)–(18), it has an important role in this paper as the basis for deriving the diffusion model. In particular, it will make clear how the singular control formulation arises. To present it, we need some notation.

Denote by $\mathcal{C}$ the set of all cycles that are subgraphs of $\mathcal{G}$, for which exactly one edge is a nonbasic activity. We call an element $c \in \mathcal{C}$ a *simple cycle* [see Figure 2(a)]. Lemma A.1 in the Appendix shows, using (12), that every nonbasic activity belongs to a simple cycle (as an edge). Consequently, there is a one-to-one correspondence between $\mathcal{E}_{nb}$ and $\mathcal{C}$, which we denote by $\sigma$. More precisely,

(20) $\quad \sigma(i,j) = c$ whenever $(i,j) \in \mathcal{E}_{nb}$ and $c$ is the simple cycle through $(i,j)$.

With an abuse of notation, we will write $(i,j) \in c$ when we mean that an (not necessarily nonbasic) activity $(i,j) \in \mathcal{E}$ belongs to the edge set of the graph $c$.

Next, we associate directions with the edges of simple cycles. Let $c$ be a simple cycle with vertices $i_0, j_0, i_1, j_1, \ldots, i_k, j_k$, where for $0 \leq l \leq k$, $i_l \in \mathcal{I}$ and $j_l \in \mathcal{J}$, and edges $(i_0, j_0) \in \mathcal{E}_{nb}$ and $(j_0, i_1), \ldots, (i_k, j_k), (j_k, i_0) \in \mathcal{E}_b$. The direction that we associate with the nonbasic element $(i_0, j_0)$ is $i_0 \to j_0$ (in words, from $i_0$ to $j_0$). The direction of the other edges, *when considered as edges of $c$*, is consistent with that of the nonbasic element, namely, $i_0 \to j_0 \to i_1 \to j_1 \to \cdots \to j_k \to i_0$. Note that an edge (corresponding to a basic activity) may have different directions when considered as an edge of different simple cycles. We signify the directions along the simple cycles by $s(c, i, j)$, defined, for all $c \in \mathcal{C}$ and $(i,j) \in c$, as

$$s(c,i,j) = \begin{cases} -1, & \text{if } (i,j), \text{ considered as an edge of } c, \text{ is directed from } i \text{ to } j, \\ 1, & \text{if } (i,j), \text{ considered as an edge of } c, \text{ is directed from } j \text{ to } i. \end{cases}$$

We will denote

(21) $\quad m_{i,c}^n = \sum_{j:(i,j) \in c} s(c,i,j) \mu_{ij}^n, \qquad i \in \mathcal{I}, \qquad m_c^n = (m_{i,c}^n, i \in \mathcal{I}).$

Next, consider the system of equations in $\psi$:

(22) $\quad \begin{cases} \sum_{j \in \mathcal{J}} \psi_{ij} = a_i, & i \in \mathcal{I}, \\ \sum_{i \in \mathcal{I}} \psi_{ij} = b_j, & j \in \mathcal{J}, \\ \psi_{ij} = 0, & (i,j) \in \mathcal{E}_{nb}. \end{cases}$

It is known that (22) has a unique solution $\psi$ for every $a, b$ satisfying $\sum_i a_i = \sum_j b_j$ (see [1], Proposition A.2). With

$$D_G = \left\{ (a,b) \in \mathbb{R}^I \times \mathbb{R}^J : \sum_{i \in \mathcal{I}} a_i = \sum_{j \in \mathcal{J}} b_j \right\}, \tag{23}$$

denote by $G : D_G \to \mathbb{R}^{IJ}$ the solution map, namely,

$$\psi_{ij} = G_{ij}(a,b), \qquad (i,j) \in \mathcal{E}. \tag{24}$$

The function $G$ is linear and so Lipschitz (a fact that will be used in the sequel). Denote also

$$H_i^n(a,b) = -\sum_j \mu_{ij}^n G_{ij}(a,b), \qquad i \in \mathcal{I}, a \in \mathbb{R}^I, b \in \mathbb{R}^J,$$

$$H^n = (H_i^n, i \in \mathcal{I}). \tag{25}$$

Let $\ell_i^n = \hat{\lambda}_i^n - \sum_{j \in \mathcal{J}} \hat{\mu}_{ij}^n \psi_{ij}^*$ and set

$$\hat{W}_i^n(t) = \hat{A}_i^n(t) - \sum_i \hat{S}_{ij}^n\left( \int_0^t \bar{\Psi}_{ij}^n(s)\, ds \right) + \ell_i^n t. \tag{26}$$

Finally, let the quantities $\{\Psi_{ij}^n\}$, as $(i,j)$ ranges over the nonbasic activities, be labeled by simple cycles, namely, define for every $c \in \mathcal{C}$, $\Psi_c = \Psi_{ij}$, where $(i,j) = \sigma^{-1}(c)$. Let a diffusion-level version of these processes be defined as $\hat{\Psi}_c^n = n^{-1/2}\Psi_c^n$.

THEOREM 2.1. *Let $X^n, Y^n, Z^n, \Psi^n$ satisfy (2)–(6) and let $\hat{X}^n, \hat{Y}^n, \hat{Z}^n, \hat{\Psi}^n$ be defined by (16)–(18). Then the following relations hold for all $t \geq 0$:*

$$\hat{X}^n(t) = \hat{X}^{0,n} + \hat{W}^n(t) + \int_0^t H^n(\hat{X}^n(s) - \hat{Y}^n(s), \hat{N}^n - \hat{Z}^n(s))\, ds$$

$$+ \sum_{c \in \mathcal{C}} m_c^n \int_0^t \hat{\Psi}_c^n(s)\, ds, \tag{27}$$

$$\hat{Y}_i^n(t) \geq 0, \qquad i \in \mathcal{I}, \qquad \hat{Z}_j^n(t) \geq 0, j \in \mathcal{J},$$

$$\sum_{i \in \mathcal{I}} [\hat{X}_i^n(t) - \hat{Y}_i^n(t)] = \sum_{j \in \mathcal{J}} [\hat{N}_j^n - \hat{Z}_j^n(t)], \tag{28}$$

$$\hat{\Psi}_c^n(t) \geq 0, \qquad c \in \mathcal{C}. \tag{29}$$

See the Appendix for a proof. We now explain how a diffusion model is derived from the above result. First, note that the limit, as $n \to \infty$, in the definition (21) of $m_c^n$ exists, and is equal to the expression obtained from (21) by replacing $\mu_{ij}^n$ by $\mu_{ij}$. We denote the limit by $m_c$. The vectors $m_c$ will be referred to as *directions of control* because of the role they will play in the singular control term. In a similar fashion, we denote by $H$ the limit of $H^n$ as $n \to \infty$, or, equivalently, as the expression obtained in (25) by replacing $\mu_{ij}^n$ by $\mu_{ij}$. Next, we take a formal limit in (26). We let $\ell = (\ell_i, i \in \mathcal{I})$, where $\ell_i = \lim_{n\to\infty} \ell_i^n = \hat{\lambda}_i - \sum_{j \in \mathcal{J}} \hat{\mu}_{ij} \psi_{ij}^*$ and let $W$ denote a Brownian motion for which the mean and the covariance matrix of $W(1)$ are given by $\ell$ and, respectively, $\Sigma := \operatorname{diag}(\lambda_i C_{U,i}^2 + \lambda_i)$ [for short, an $(\ell, \Sigma)$-Brownian motion]. The expression above for the covariance is obtained by calculating the covariance matrix of $\hat{W}^n(1)$ upon formally replacing the fluid-level quantities $\bar{\Psi}_{ij}^n(s)$ in (26) by the quantities $\psi_{ij}^*$ from the fluid model, and finally taking limits as $n \to \infty$.

Next, equation (29) imposes a constraint on $\hat{\Psi}_c^n$ in the form of a lower bound at zero. It is also not hard to see that an upper bound on $\hat{\Psi}_c^n$ follows from (5), (14) and (18). However, such a bound is of the order of $n^{1/2}$, and as $n$ grows without bound, it becomes irrelevant. It therefore makes sense that, in the proposed diffusion model, each integral in the last term of (27) is replaced by a process $\eta_c$ that is required to have increasing sample paths, and no additional limitation.

We thus obtain a diffusion model involving processes $X(t)$, $Y(t)$, $Z(t)$ taking values in $\mathbb{R}_+^I$, $\mathbb{R}_+^I$ and, respectively, $\mathbb{R}_+^J$ as well as $\mathbb{R}_+$-valued processes $\eta_c$, $c \in \mathcal{C}$. Equations (30)–(32) below, that describe the diffusion model, are analogous to equations (27)–(29), respectively:

$$(30) \quad X(t) = x + W(t) + \int_0^t H(X(s) - Y(s), -Z(s))\, ds + \sum_{c \in \mathcal{C}} m_c \eta_c(t)$$

$$(31) \quad \begin{aligned} Y_i(t) &\geq 0, \quad i \in \mathcal{I}, \quad Z_j(t) \geq 0, \quad j \in \mathcal{J}, \\ \sum_{i \in \mathcal{I}} Y_i(t) &= \sum_{i \in \mathcal{I}} X_i(t) + \sum_{j \in \mathcal{J}} Z_j(t), \quad t \geq 0. \end{aligned}$$

(32)    For every $c \in \mathcal{C}$, $\eta_c$ is nondecreasing and $\eta_c(0) \geq 0$.

We will consider $Y, Z$ and $\{\eta_c\}$ as control processes, and regard (31), (32) as constraints that they must satisfy. To define precisely what controls will be regarded as admissible, recall that $W$ is an $(\ell, \Sigma)$-Brownian motion defined on $(\Omega, \mathcal{F}, P)$. Let $(\mathcal{F}_t)$ be a right-continuous filtration of sub-sigma-fields of $\mathcal{F}$ such that $(W(t) - \ell t, \mathcal{F}_t : t \geq 0)$ is a martingale. Let a deterministic initial condition $x \in \mathbb{R}^I$ be given. A triplet $(Y, Z, \eta)$ is said to be an *admissible* *control* and $X$ a corresponding *controlled process* if $Y$, $Z$ and $\eta$ are processes with cadlag sample paths, $Y$ and $Z$ are $(\mathcal{F}_t)$-progressively measurable, $\eta_c$ is $(\mathcal{F}_t)$-adapted, and (30)–(32) hold for all $t \in [0, \infty)$, $P$-a.s.

A principal hypothesis of all the results below is what we refer to as the *null controllability condition*:

(33)             There exists $c \in \mathcal{C}$ such that $e \cdot m_c < 0$.

Note that the cone $\mathbb{C}$ referred to in the Introduction, in which the increments of the $I$-dimensional process $\sum_c m_c \eta_c(t)$ of (30) take values, is simply the cone generated by $\{m_c : c \in \mathcal{C}\}$. Also, condition (1) of the Introduction is given explicitly by (33). The following result is proved in Section 3.

THEOREM 2.2. *Assume (33) holds. Then there exists an admissible control $(Y, Z, \eta)$ under which $e \cdot X(t) \leq 0$ and $Y(t) = 0$ on $[0, \infty)$, $P$-a.s.*

The main results of this paper establish the validity of statements about the queueing model, which are analogous to Theorem 2.2 in an asymptotic sense. The first is concerned with preemptive scheduling (see Section 3 for a proof).

THEOREM 2.3. *Assume (33) holds. Then there exist P-SCPs under which, for every $\varepsilon$ and $T$ satisfying $0 < \varepsilon < T < \infty$,*

$$\lim_{n \to \infty} P(Y^n(t) = 0 \text{ for all } t \in [\varepsilon, T]) = 1.$$

The treatment of nonpreemptive scheduling is more involved. In order to keep the notation simple, we have focused in this paper on the most simple case where one can expect a null-controllability result: The case $I = J = 2$. Clearly, in this case there is at most one simple cycle. It is expected that the general case can be treated with similar ideas. The result below is proved in Section 4.

THEOREM 2.4. *Assume $I = J = 2$ and let (33) hold. Then there exist N-SCPs under which, for every $\varepsilon$ and $T$ satisfying $0 < \varepsilon < T < \infty$,*

$$\lim_{n \to \infty} P(Y^n(t) = 0 \text{ for all } t \in [\varepsilon, T]) = 1.$$

REMARK. As can be seen in Sections 3 and 4, Theorems 2.3 and 2.4 hold with $\varepsilon = 0$ in the case that the initial condition $x$ satisfies $e \cdot x < -1$ (or even $e \cdot x < 0$, with an obvious modification of the proofs).

2.4. *Discussion.* Let us consider some numerical examples. Consider first a system with structure as depicted in Figure 1(a). Assume $\nu_j = 1$ for $j = 1, 2, 3$, and

$$\lambda = \begin{pmatrix} 8 \\ 4 \end{pmatrix}, \qquad \bar{\mu} = \mu = \begin{pmatrix} 3 & 10 & 1 \\ 1 & 4 & 2 \end{pmatrix},$$

(where in this subsection we abuse notation and label $j$ by $1, \ldots, J$). One checks that the heavy traffic condition holds, and

$$\xi^* = \begin{pmatrix} 1 & 0.5 & 0 \\ 0 & 0.5 & 1 \end{pmatrix}, \qquad m_1 = \begin{pmatrix} -7 \\ 3 \end{pmatrix}, \qquad m_2 = \begin{pmatrix} 9 \\ -2 \end{pmatrix}.$$

Above, $m_1$ and $m_2$ are the directions of control corresponding to the two nonbasic activities $(2, 1)$ and $(1, 3)$. In fact, the graph that appears in Figure 1(b) precisely describes $\mathcal{G}$ and $\mathcal{G}_b$ in the current example. Clearly, the null-controllability condition holds, since $e \cdot m_1 < 0$. The simple cycles and directions of control are depicted in Figure 2(a) and (b). To demonstrate the geometric aspect, we refer to Figure 2(c), where the null domain $\mathbb{X}_-$ is shown along with a vector field defined on its boundary assuming the constant value $m_1$. Under appropriate assumptions, one can construct a diffusion in $\mathbb{R}^2$ with a boundary term according to this vector field, that will be constrained to $\mathbb{X}_-$. In contrast, $m_2$ cannot be used to constrain the diffusion to the same domain. In general, the diffusion can be constrained to $\mathbb{X}_-$, provided that at least one of the vectors $m_c$ satisfies $e \cdot m_c < 0$. This is the source for condition (33).

We next consider examples with two classes and two stations:

$$\mu = \bar{\mu} = \begin{pmatrix} 8 & 10 \\ 3 & 6 \end{pmatrix}, \qquad \lambda = \begin{pmatrix} 13 \\ 3 \end{pmatrix},$$

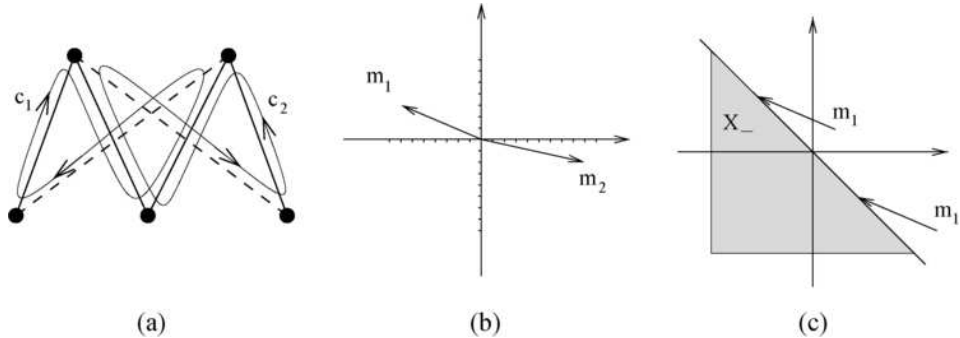

FIG. 2. *(a) A graph as in Figure* 1(b) *with two simple cycles. Dashed lines represent nonbasic activities. (b) A possible set of directions of control corresponding to the two cycles. (c) The direction $m_1$ may be used to constrain the diffusion model to the null domain $\mathbb{X}_-$.*

(34)
$$\xi^* = \begin{pmatrix} 1 & 0.5 \\ 0 & 0.5 \end{pmatrix}, \qquad m = \begin{pmatrix} -2 \\ 3 \end{pmatrix},$$

$$\mu = \bar{\mu} = \begin{pmatrix} 4 & 7 \\ 2 & 4 \end{pmatrix}, \qquad \lambda = \begin{pmatrix} 7.5 \\ 2 \end{pmatrix},$$

(35)
$$\xi^* = \begin{pmatrix} 1 & 0.5 \\ 0 & 0.5 \end{pmatrix}, \qquad m = \begin{pmatrix} -3 \\ 2 \end{pmatrix}.$$

Above, we have presented both the data and the solution to the linear program in each case. In both cases the heavy traffic condition holds and the activities are as depicted in Figure 3(a). Note that in these examples there is a single cycle. The null controllability condition does not hold in the first example, and it does hold in the second. In fact, one can write the null controllability condition (33) for the above examples in a straightforward way as

(36) $$\mu_{11} + \mu_{22} < \mu_{12} + \mu_{21}.$$

It is instructive to note that there are values of $\mu$ for which (36) does not express the null controllability condition (33). Consider the following example:

$$\mu = \bar{\mu} = \begin{pmatrix} 3 & 7 \\ 6 & 11 \end{pmatrix}, \qquad \lambda = \begin{pmatrix} 3.5 \\ 11.5 \end{pmatrix},$$

$$\xi^* = \begin{pmatrix} 0 & 0.5 \\ 1 & 0.5 \end{pmatrix}, \qquad m = \begin{pmatrix} 4 \\ -5 \end{pmatrix}.$$

Note that the values of both sides of the inequality (36) are the same as in example (34). Still, as can be verified by calculating $e \cdot m$, the null controllability condition, that did not hold in example (34), does hold in this example. The reason is that the structure of the graph has changed and it now corresponds to Figure 3(b). In particular, the direction of the cycle is reversed, and condition (33) is equivalent to (36) with a reversed inequality. We can see that the null controllability condition depends on the values of $\mu$, as well as the direction of the cycle, which in turn is determined by the fluid model parameters $\lambda$ and $\bar{\mu}$.

We now come back to the point referred to in the Introduction regarding the unusual heavy traffic behavior. A particular consequence of our main results can be stated as follows. One can find policies (preemptive and nonpreemptive) under which, for every $t > 0$,

(37) $$\lim_{n \to \infty} P(e \cdot Y^n(t) = 0) = 1.$$

What we referred to as unusual is that a critically loaded system shows a behavior that is typical to underloaded systems: the possibility to maintain empty queues with probability approaching one. We would like to make precise the statement that the system under study is critically loaded. One aspect of this is that the underlying fluid model is critically loaded, in the sense that the linear program (10) is solved with $\rho^* = 1$, by assumption. More significant, however, is a statement that can be made regarding the probabilistic model. Namely, one can show that in a probabilistic model associated with an overloaded fluid model, queues inevitably build up.

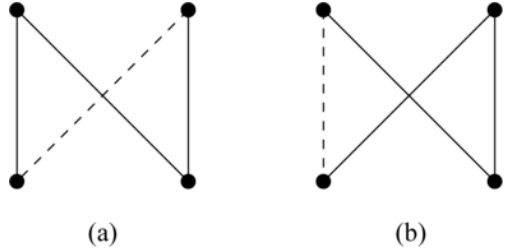

Fig. 3. *Two examples with cycles in opposite directions (dashed lines represent nonbasic activities).*

To this end, let $\lambda'_i$, $i \in \mathcal{I}$ be constants satisfying $\lambda'_i \geq \lambda_i$ for all $i$ and $\lambda'_i > \lambda_i$ for at least one $i \in \mathcal{I}$. The parameters $(\lambda, \bar\mu)$ lead to a fluid model that is critically loaded. In the same sense, the pair $(\lambda', \bar\mu)$ correspond to an overloaded fluid model. We consider a sequence of processes $A^n_{OL}$ defined analogously to $A^n$, but with $\lambda'^n$ replacing the parameters $\lambda^n$, where $\lambda'^n$ is a sequence satisfying $n^{-1}\lambda'^n \to \lambda'$ [compare with (8)]. Let $X^n_{OL}$ stand for the processes $X^n$ obtained by replacing throughout our probabilistic model $A^n$ by $A^n_{OL}$. Define analogously all other processes involved, for example, $Y^n_{OL}$ in place of $Y^n$. As the following result shows, the model thus obtained is indeed overloaded in the sense that queues (in fact, large queues) necessarily build up. The result shows a sharp contrast with (37).

PROPOSITION 2.1. *There exist constants $C_1$ and $C_2 > 0$ depending only on $(\lambda, \lambda', \bar\mu)$ (and not depending on $n$ or $t$) such that, under any policy, for every $t$,*

(38) $$\lim_{n \to \infty} P(e \cdot Y^n_{OL}(t) \geq (-C_1 + C_2 t)n) = 1.$$

A proof is sketched in the Appendix.

**3. Diffusion model and queueing model in the preemptive case.** In this section we prove Theorems 2.2 and 2.3.

PROOF OF THEOREM 2.2. Let $c_0$ be such that $e \cdot m_{c_0} < 0$. Consider the domain

$$D_\alpha = \{\xi \in \mathbb{R}^I : e \cdot \xi < -\alpha\}$$

for some fixed $\alpha \geq 0$. Fix some $j_0 \in \mathcal{J}$. We will construct a control with the following properties:

(39)
$$Y(t) = 0, \qquad Z_j(t) = 0, \qquad \text{for all } j \neq j_0,$$
$$\eta_c(t) = 0 \qquad \text{for all } c \neq c_0, \quad t \geq 0.$$

The process $Z_{j_0}$ will satisfy

(40)
$$Z_{j_0}(t) = -e \cdot X(t), \qquad t \geq 0.$$

As a result, (31) will be satisfied. Finally, the process $\eta_{c_0}$ will serve as a constraining term of a reflected diffusion on the domain $D_\alpha$ with reflection field identically equal to $m_{c_0}$ on the boundary $\partial D_\alpha$. We therefore consider equation (30) in the special form

(41)
$$X(t) = x + W(t) + \int_0^t \widetilde{H}(X(s))\,ds + m_{c_0}\eta_{c_0}(t),$$

where $\widetilde{H}(\xi) = H(\xi, (e \cdot \xi)e_{j_0})$. Note that $\widetilde{H}$ is linear. Consider first the case $x \in \overline{D_\alpha}$. The result of [10] regarding existence of strong solutions to stochastic differential equations with oblique reflection, stated for a bounded domain, implies, using a standard localization argument, the existence of a pair $(X, \eta_{c_0})$ with the following properties. The process $X$ is progressively measurable, $\eta_{c_0}$ is continuous nondecreasing, adapted, with values in $\mathbb{R}_+$, equation (41) holds for all $t \geq 0$, a.s. and $X(t) \in \overline{D_\alpha}$, $t \geq 0$, a.s. In particular, we have constructed a process $X$ with the property $e \cdot X(t) \leq -\alpha$, $t \geq 0$, a.s. Letting now $Y$, $Z$ and $\eta_c$, $c \neq c_0$ be defined via (39) and (40), we have constructed a triplet $(Y, Z, \eta)$ that is an admissible control, and have shown that $X$ and $Y$ satisfy the conclusion of the theorem.

In the case where $e \cdot x > -\alpha$, clearly, $x + \beta m_{c_0} \in \overline{D_\alpha}$ for $\beta$ large. Fix any $\beta$ as above, and set $y = x + \beta m_{c_0} \in \overline{D_\alpha}$. Denote the control and the controlled process corresponding to starting from $y$ as $(\widetilde{Y}, \widetilde{Z}, \widetilde{\eta})$ and, respectively, $\widetilde{X}$. Now set $Y = \widetilde{Y}$, $Z = \widetilde{Z}$, $\eta_{c_0}(t) = \beta + \widetilde{\eta}_{c_0}(t)$ and $X(t) = \widetilde{X}(t)$. Then $X(0) = y$, and clearly, (41) still holds for all $t \geq 0$. As a result, $(Y, Z, \eta)$ is an admissible control and the conclusion of the theorem holds. $\square$

We remark that in the proof above we can simply take $\alpha = 0$. Our results regarding asymptotic null controllability will be inspired by a similar construction, but $\alpha$ will be taken to be positive.

Recall that, by definition, a P-SCP is a map that determines $\Psi^n(t)$ for a given value of $X^n(t)$ in such a way that (2)–(6) hold. The following lemma shows that, under suitable conditions, we can determine the value of $\Psi^n$ by first selecting values for $Y^n(t)$, $Z^n(t)$ and $\{\Psi_c^n(t), c \in \mathcal{C}\}$.

LEMMA 3.1. *There is a constant $a_0 > 0$ such that the following holds for all $n$ large. Suppose that, for some $t$, the following relations hold: $X^n(t) \in \mathbb{Z}_+^I$, $\|X^n(t) - nx^*\| \leq a_0 n$, $Y^n(t) \in \mathbb{Z}_+^I$, $Z^n(t) \in \mathbb{Z}_+^J$,*

$$e \cdot Y^n(t) + e \cdot N^n = e \cdot X^n(t) + e \cdot Z^n(t),$$

$$[e \cdot Y^n(t)] \vee [e \cdot Z^n(t)] \leq a_0 In + 1,$$

$\Psi_c^n(t) \in \mathbb{Z}_+$, $c \in \mathcal{C}$, *and*

$$\Psi_c^n(t) \leq a_0 n, \qquad c \in \mathcal{C}.$$

*Set*

(42) $\quad \Psi_{ij}^n(t) = G_{ij}(X^n(t) - Y^n(t), N^n - Z^n(t)) - \sum_{c \in \mathcal{C} : (i,j) \in c} s(c,i,j) \Psi_c^n(t).$

*Then the quantities*

$$\{X_i^n(t), i \in \mathcal{I}\}, \quad \{Y_i^n(t), i \in \mathcal{I}\}, \quad \{Z_j^n(t), j \in \mathcal{J}\}, \quad \{\Psi_{ij}^n(t), (i,j) \in \mathcal{I} \times \mathcal{J}\}$$

*satisfy (2) and (4)–(6).*

The proof of this lemma appears in the Appendix.

PROOF OF THEOREM 2.3. We begin by describing the policy. By Definition 2.1, we need to describe a map that determines $\Psi^n(t)$ for a given value of $X^n(t)$. Given the one-to-one relations (16) and (18), this is equivalent to describing $\hat{\Psi}^n(t)$ for a given value of $\hat{X}^n(t)$. Let $a_0$ be the constant from Lemma 3.1. If $\|\hat{X}^n(t)\| > a_0 n^{1/2}$, we assume that $\hat{\Psi}^n(t)$ is determined as the image of $\hat{X}^n(t)$ under some fixed map $f^n$, in a way that is consistent with Definition 2.1 but otherwise arbitrary (there will be no need to describe $f^n$ more precisely). Focusing below on the case $\|\hat{X}^n(t)\| \leq a_0 n^{1/2}$, we fix a sequence $K_n$, $n \in \mathbb{N}$, so that $n^{1/2} K_n \in \mathbb{Z}_+$ for all $n$ and

(43) $\qquad K_n \to \infty, \qquad n^{-1/2} K_n \to 0 \qquad \text{as } n \to \infty.$

Denote

$$D_1 = \{\xi \in \mathbb{R}^I : e \cdot \xi < -1\}.$$

Fix throughout a simple cycle $c_0$ for which $e \cdot m_{c_0} < 0$. Also fix throughout $i_0 \in \mathcal{I}$ and $j_0 \in \mathcal{J}$. The proposed policy sets

(44) $\qquad \hat{\Psi}_{c_0}^n(t) = \begin{cases} 0, & \hat{X}^n(t) \in D_1, \\ K_n, & \hat{X}^n(t) \in D_1^c, \end{cases} \qquad t \geq 0,$

and $\hat{\Psi}_c^n(t) = 0$ for all $c \neq c_0$, $t \geq 0$. It also sets $\hat{Y}_i^n(t) = 0$ for all $i \neq i_0$, $\hat{Z}_j^n(t) = 0$ for all $j \neq j_0$ and

$$
\begin{aligned}
\hat{Y}_{i_0}^n(t) &= (e \cdot \hat{X}^n(t) - e \cdot \hat{N}^n)^+, \\
\hat{Z}_{j_0}^n(t) &= (e \cdot \hat{X}^n(t) - e \cdot \hat{N}^n)^-, \qquad t \geq 0.
\end{aligned}
\tag{45}
$$

By (43) and (44), $\Psi_{c_0}^n(t) \leq a_0 n$ for all $n$ large. By Lemma 3.1, $\Psi^n(t)$ are well defined given $X^n(t)$, and (2), (4)–(6) hold. This completes the description of the P-SCP. Clearly, this description, along with equation (3), uniquely defines the processes $\Psi^n$ and $X^n$ given the initial conditions and primitive processes. Although the description applies for any initial condition $\hat{X}^{0,n}$, the treatment with be slightly different for $\hat{X}^{0,n}$ in $D_1$ and outside.

In what follows let $0 < \varepsilon < T < \infty$ be fixed. Let

$$
\bar{\tau}^n = \inf\{t : \|\hat{X}^n(t)\| > a_0 n^{1/2}\}.
\tag{46}
$$

Let also

$$
\widetilde{H}^n(\xi) = H^n(\xi - (e \cdot \xi - e \cdot \hat{N}^n)^+ e_{i_0}, -(e \cdot \xi - e \cdot \hat{N}^n)^- e_{j_0}).
$$

For $t \leq \bar{\tau}^n$, one has by (27)

$$
\hat{X}^n(t) = \hat{X}^{0,n} + \hat{W}^n(t) + \int_0^t \widetilde{H}^n(\hat{X}^n(s)) \, ds + m_{c_0}^n \int_0^t \hat{\Psi}_{c_0}^n(s) \, ds.
\tag{47}
$$

Note that $\widetilde{H}^n$ satisfy

$$
\|\widetilde{H}^n(\xi)\| \leq C_H(\|\xi\| + \|\hat{N}^n\|), \qquad \xi \in \mathbb{R}^I, n \in \mathbb{N},
\tag{48}
$$

where $C_H$ is a constant independent of $n$. Denote $C_e^n = -e \cdot m_{c_0}^n$. Note that, by assumption, $C_e^n \to -e \cdot m_{c_0} > 0$ as $n \to \infty$. It is assumed in what follows that $n$ is so large that

$$
C_e^n > |e \cdot m_{c_0}|/2 =: 2C_0.
\tag{49}
$$

Denote $\hat{X}_e^n = e \cdot \hat{X}^n$ and, similarly, $\hat{W}_e^n = e \cdot \hat{W}^n$, $\hat{X}_e^{0,n} = e \cdot \hat{X}^{0,n}$ and $\widetilde{H}_e^n = e \cdot \widetilde{H}^n$. Then by (44) and (47), we have

$$
\begin{aligned}
\hat{X}_e^n(t) &= \hat{X}_e^{0,n} + \hat{W}_e^n(t) + \int_0^t \widetilde{H}_e^n(\hat{X}^n(s)) \, ds \\
&\quad - C_e^n K_n \int_0^t \mathbf{1}_{\hat{X}_n^e(s) \geq -1} \, ds, \qquad t \leq \bar{\tau}^n.
\end{aligned}
\tag{50}
$$

The rest of our argument is divided into four steps.

*Step* 1. We first show that there exists a deterministic constant $C_1$ independent of $x$, $n$ and $K_n$ such that

(51) $\qquad \|\hat{X}^n\|_T^* \leq C_1(1 + \|\hat{X}^{0,n}\| + \|\hat{W}^n\|_T^*) \qquad$ on the event $\bar{\tau}^n \geq T$.

To this end, denote $A_n = [4\|e\|(\|\hat{X}^{0,n}\| + \|\hat{W}^n\|_T^*)] \vee 1$ and let

$$\tau_1^n = \inf\{t \in [0,T] : \hat{X}_e^n(t) \leq -A_n - 1\}.$$

Note that $\tau_1^n > 0$. Since by assumption $\bar{\tau}^n \geq T$, (47) and (50) are valid for $t \leq T$. By (48) and (50), noting that $\hat{X}_e^n(t) \geq -A_n - 1$ for $t \leq \tau_1^n$, we have

$$K_n C_e^n \int_0^t \mathbf{1}_{\hat{X}_e^n(s) \geq -1}\, ds \leq A_n + 1 + \beta_n + \|e\|\|\hat{X}^{0,n}\| + \|e\|\|\hat{W}^n\|_T^*$$
$$+ \|e\|C_H \int_0^t \|\hat{X}^n(s)\|\, ds, \qquad t \leq \tau_1^n \wedge T,$$

where we denote $\beta_n = \|e\|C_H T\|\hat{N}^n\|$. Note that $\beta_n \to 0$. Hence, by (47),

$$\|\hat{X}^n(t)\| \leq C_2(A_n + 1) + C_2 \int_0^t \|\hat{X}^n(s)\|\, ds, \qquad t \leq \tau_1^n \wedge T,$$

where $C_2$ does not depend on $\hat{X}^{0,n}$, $\hat{W}^n$, $K_n$ or $n$. By Gronwall's lemma,

(52) $\qquad\qquad\qquad \|\hat{X}^n\|_{\tau_1^n \wedge T}^* \leq C_2(A_n + 1)e^{C_2 T}.$

Since on the event $\tau_1^n \geq T$ the property (51) follows from (52), we consider in what follows only the case $\tau_1^n < T$. If $\hat{X}_e^n(t) < -1$ for all $t \in [0, \tau_1^n]$, let $\tau_2^n = 0$; otherwise let

$$\tau_2^n = \sup\{t \in [0, \tau_1^n] : \hat{X}_e^n(t) \geq -1\}.$$

Note that $\hat{X}_e^n(\tau_2^n) \geq -1 - \|e\|\|\hat{X}^{0,n}\|$. Hence, by (50),

$$-A_n + \|e\|\|\hat{X}^{0,n}\| \geq \hat{X}_e^n(\tau_1^n) - \hat{X}_e^n(\tau_2^n)$$
$$= \hat{W}_e^n(\tau_1^n) - \hat{W}_e^n(\tau_2^n) + \int_{\tau_2^n}^{\tau_1^n} \widetilde{H}_e^n(\hat{X}^n(s))\, ds.$$

As a result,

$$A_n/2 \leq \left|\int_{\tau_2^n}^{\tau_1^n} \widetilde{H}_e^n(\hat{X}^n(s))\, ds\right| \leq \|e\|C_H \tau_1^n \|\hat{X}^n\|_{\tau_1^n}^* + \beta_n.$$

This and (52) show $A_n/2 - \beta_n \leq \|e\|C_H \tau_1^n C_2(A_n + 1)e^{C_2 T}$. Since $A_n \geq 1$, $A_n/(A_n + 1) \geq 1/2$. As a result, there exists a deterministic constant $\delta > 0$ not depending on $\hat{X}^{0,n}, \hat{W}_n, K_n, n$ such that, provided that $n$ is large enough,

(53) $\qquad\qquad\qquad \tau_1^n \geq \delta \qquad$ on the event $\bar{\tau}^n \geq T$.

Combining (52) and (53), we have that $\|\hat{X}^n\|_\delta^* \leq C_3(1 + \|\hat{X}^{0,n}\| + \|\hat{W}^n\|_T^*)$, where $C_3$ does not depend on $\hat{X}^{0,n}, \hat{W}^n, K_n, n$. In a similar fashion we have, for $i = 1, 2, \ldots$,

$$\|\hat{X}^n\|_{i\delta}^* \leq C_3(1 + \|\hat{X}^n\|_{(i-1)\delta}^* + \|\hat{W}^n\|^*).$$

Assuming without loss that $C_3 \geq 1$, the above implies

$$\|\hat{X}^n\|_T^* \leq (2C_3)^{T/\delta+1}(1 + \|\hat{X}^{0,n}\| + \|\hat{W}^n\|_T^*)$$

and proves (51).

*Step* 2. Recall that $W$ denotes an $(\ell, \Sigma)$-Brownian motion. We next show that

(54) $$\lim_{n \to \infty} P(\bar{\tau}^n \leq T) = 0$$

and

(55) $$\hat{W}^n \Rightarrow W \qquad \text{on } [0, T].$$

Let $A_i$, $i \in \mathcal{I}$, and $S_{ij}$, $(i,j) \in \mathcal{I} \times \mathcal{J}$ be mutually independent Brownian motions with mean zero and variances given by $EA_i^2(1) = \lambda_i C_{U,i}^2$, $ES_{ij}^2(1) = \mu_{ij}$. By Theorem 14.6 of [4],

(56) $$(\hat{A}^n, \hat{S}^n) \Rightarrow (A, S) \qquad \text{locally on compacts.}$$

By (5) and (6), $\Psi_{ij}^n(t) \leq N_j^n$ for all $i, j$ and $t$. Thus, by (8) and (26), $\|\hat{W}^n\|_T^* \leq \|\hat{A}^n\|_T^* + \|\hat{S}^n\|_{C_4T}^* + \|\ell^n\|T$ for a suitable constant $C_4$. Hence, $n^{-1/2}\|\hat{W}^n\|_T^*$ converges to zero in probability. By (51), $n^{-1/2}\|\hat{X}^n\|_{T\wedge\bar{\tau}^n}^*$ converges to zero in probability. Using the definition (46), this establishes (54). In turn, this shows that $n^{-1/2}\|\hat{X}^n\|_T^*$ converges to zero in probability. By (45), so do $n^{-1/2}\|\hat{Y}^n\|_T^*$ and $n^{-1/2}\|\hat{Z}^n\|_T^*$. Note that $G(x^*, \nu) = \psi^*$, as follows from (11) and (24). Using (42), linearity of the map $G$, (43) and (16)–(18), we thus obtain, for a suitable constant $C_5$,

$$\|\bar{\Psi}^n(t) - \psi^*\|_T^* \leq n^{-1}\|G(X^n - nx^* - Y^n, N^n - n\nu - Z^n)\|_T^* + n^{-1/2}K_n$$
(57) $$\leq C_5 n^{-1/2}(\|\hat{X}^n\|_T^* + \|\hat{Y}^n\|_T^* + \|\hat{N}^n\| + \|\hat{Z}^n\|_T^*) + n^{-1/2}K_n$$
$$\to 0 \qquad \text{in probability.}$$

Combining (26), (56) and (57), the claim (55) follows from the lemma on page 151 of [4] regarding random change of time.

*Step* 3. Recall (19). We prove the theorem in the case $x \in D_1$. In this case, for all $n$ large, $\hat{X}^{0,n} \in D_1$. Denote

$$\tau^n = \inf\{s \in [0,T] : \hat{X}^n_e(s) \geq -1/2\}.$$

Denote also $\Omega^n = \{\bar{\tau}^n > T\}$. By (54), $\lim_n P(\Omega^n) = 1$. If we show

(58)
$$\lim_{n \to \infty} P(\{\tau^n \leq T\} \cap \Omega^n) = 0,$$

it would follow that $P(\hat{X}^n_e(t) \leq -1/2$ for all $t \in [0,T]) \to 1$ and, in turn, by (17) and (45), that $P(Y^n(t) = 0$ for all $t \in [0,T]) \to 1$, as $n \to \infty$. Hence, in order to prove the theorem, it suffices to prove (58).

To this end, note that the jumps of the process $\hat{W}^n$ are bounded by $n^{-1/2}$ and write

(59)
$$\begin{aligned}P(\{\tau^n \leq T\} \cap \Omega^n) \leq P(\{&\text{there exist } 0 \leq s < t \leq T \text{ such that}\\ &-1 \leq \hat{X}^n_e(\theta) \leq -1/2 \text{ for all } \theta \in [s,t],\\ &\hat{X}^n_e(s) \leq -7/8 \text{ and } \hat{X}^n_e(t) \geq -5/8\} \cap \Omega^n).\end{aligned}$$

Under the event indicated immediately above, on the (random) interval $[s,t]$, $\hat{X}^n_e \geq -1$, and thus, by (50),

(60)
$$\begin{aligned}\tfrac{1}{4} &\leq \hat{X}^n_e(t) - \hat{X}^n_e(s)\\ &= \int_s^t \widetilde{H}^n_e(\hat{X}^n(\theta))\,d\theta + \hat{W}^n_e(t) - \hat{W}^n_e(s) - C^n_e K_n(t-s).\end{aligned}$$

Moreover, by (48) and (51),

$$\int_s^t \widetilde{H}^n_e(\hat{X}^n(\theta))\,d\theta \leq C_H \|e\| C_1 (t-s)[r + \|W_n\|^*_T] + \beta_n,$$

where $r = 1 + \sup_n \|\hat{X}^{0,n}\|$. Since $K_n \to \infty$, it follows that there is a deterministic $n_0$ such that, for all $n \geq n_0$,

(61)
$$\begin{aligned}\Delta_n(s,t) &:= \hat{W}^n_e(t) - \hat{W}^n_e(s) + C_H\|e\|C_1(t-s)\|\hat{W}^n\|^*_T + \beta_n\\ &\geq \tfrac{1}{4} + [C^n_e K_n - C_H\|e\|C_1 r](t-s)\\ &\geq \tfrac{1}{4} + C_0 K_n(t-s)\\ &\geq \begin{cases} C_0 K_n^{1/2}, & t - s \geq K_n^{-1/2},\\ 1/4, & t - s < K_n^{-1/2},\end{cases}\end{aligned}$$

where $C_0$ is as in (49), and on the first inequality we used the fact that $C_H\|e\|C_1 r$ does not depend on $n$. Combining (59) and (61), we see that there are constants $C_6, C_7 > 0$ not depending on $n$ and $K_n$ such that

$$P(\{\tau^n \leq T\} \cap \Omega^n)$$
$$\leq P(\text{there exist } s < t \leq T \text{ such that } \Delta_n(s,t) \geq C_0 K_n^{1/2})$$
$$\quad + P(\text{there exist } s < t \leq T \text{ such that } t - s < K_n^{-1/2}, \Delta_n(s,t) \geq 1/4)$$
$$\leq P(\|\hat{W}_n\|_T^* \geq C_6 K_n^{1/2})$$
$$\quad + P(w_T(\hat{W}_n, K_n^{-1/2}) \geq 1/8 - \beta_n)$$
$$\quad + P(C_H \|e\| C_1 K_n^{-1/2} \|\hat{W}_n\|_T^* \geq 1/8)$$
$$\leq 2P(\|\hat{W}_n\|_T^* \geq C_7 K_n^{1/2}) + P(w_T(\hat{W}_n, K_n^{-1/2}) \geq 1/8 - \beta_n),$$

where $w_T(f, \delta)$ denotes the modulus of continuity of a function $f$ on $[0,T]$. Since by (55) $\hat{W}^n$ are tight, and the weak limit process has continuous sample paths, (58) follows (see e.g., Theorem 16.8 of [4]) and, hence, the result.

*Step* 4. Finally, we prove the theorem in the case $x \in D_1^c$. Let $\tau_3^n = \inf\{t \in [0,T] : \hat{X}_e^n < -1\}$. Then

$$P(\{\tau_3^n > \varepsilon\} \cap \Omega^n) \leq P(\{\hat{X}_e^n(\theta) \geq -1 \text{ for all } \theta \in [0,\varepsilon]\} \cap \Omega^n).$$

An argument similar to step 3 (only simpler) shows that, under the event indicated on the r.h.s. of the above display,

$$\Delta_n(0,\varepsilon) \geq -1 - r + [C_e^n K_n - C_H \|e\| C_1 r]\varepsilon,$$

and, in turn, that $P(\tau_3^n > \varepsilon)$ converges to zero as $n \to \infty$. We can now review the argument of step 3, replacing $\hat{X}^{0,n}$ by $\hat{X}^n(\tau_3^n)$ and $\tau^n$ by $\tau_4^n := \inf\{t \in [\tau_3^n, T] : \hat{X}_e^n(t) \geq -1/2\}$ (where $\tau_4^n = \infty$ on the event $\tau_3^n > T$) so as to show that $P(\tau_4^n \leq T) \to 0$ as $n \to \infty$. The only issue that is different now is that the "initial condition" $\hat{X}^n(\tau_3^n)$ is random and cannot be bounded by a deterministic constant $r$. Instead, let us define the random variable $r^n := 1 + \|\hat{X}^n(\tau_3^n)\|$ and let $\Omega_1^n := \{r^n \leq C_0(C_H \|e\| C_1)^{-1} K_n\} \cap \{\tau_3^n \leq T\}$. Coming back to (61) with $r^n$ in place of $r$, it is clear that the second inequality will hold on $\Omega_1^n$, and so the remainder of the argument of step 3 is valid once $\Omega^n$ is replaced by $\Omega^n \cap \Omega_1^n$. Since $P(\tau_3^n \leq T)$ converges to one, the relations (51), (54) and (55) imply that the random variables $r^n$ are tight and, therefore, $P(\Omega_1^n \cap \Omega^n) \to 1$ as $n \to \infty$. This establishes the theorem. □

**4. The nonpreemptive case.** In this section we treat the nonpreemptive case and prove Theorem 2.4. Recall that we only consider the case $I = J = 2$. Thus, $\mathcal{I} = \{1,2\}$ and $\mathcal{J} = \{3,4\}$. The heavy traffic and complete resource pooling conditions, which are in force, imply that the graph of basic activities $\mathcal{G}_b$ is a tree with vertex set $\mathcal{I} \cup \mathcal{J} = \{1, 2, 3, 4\}$ [cf. (12)]. It follows that there are exactly 3 activities that are basic. Had we not had a fourth activity, the graph $\mathcal{G}$ would not contain a cycle and it would not be possible to fulfill the null controllability condition (33). Thus, the hypotheses of the theorem require that there be a fourth, nonbasic activity. We have labeled the classes as "1" and "2." Although $\mathcal{J} = \{3, 4\}$, it will be more natural in the discussion that follows, and throughout the section, to refer to the stations as "station 1" and "station 2," and with an abuse of notation, regard the index set for the stations as $\{1, 2\}$. Accordingly, we have four activities, $(i, j)$, $i, j \in \{1, 2\}$, and without loss of generality, we let $(2, 1)$ be the only nonbasic activity [see Figure 3(a)]. As a result, the direction of the only simple cycle, which we denote as $c$, is the following: class $2 \to$ station $1 \to$ class $1 \to$ station $2 \to$ class 2. By (21), we get $m_c^n = (\mu_{11}^n - \mu_{12}^n, -\mu_{21}^n + \mu_{22}^n)$ and $m_c = (\mu_{11} - \mu_{12}, -\mu_{21} + \mu_{22})$. In what follows we will write $m^n$ for $m_c^n$ and $m$ for $m_c$. We also let $C_m^n = -e \cdot m^n$ and $C_m = -e \cdot m$. Note that condition (33) can be written as $C_m \equiv -e \cdot m > 0$.

We now specify the N-SCP for which the conclusions of Theorem 2.4 will be shown. According to Definition 2.2, we must specify the initial arrangement and how routing is determined upon each arrival. To this end, we need some notation. Note that by (23) and (25) there exists a constant $C_H'$ such that

$$\|H^n(a, b)\| \leq \frac{C_H'}{2I}(\|a\| + \|b\|), \qquad (a, b) \in D_G, n \in \mathbb{N}. \tag{62}$$

Let

$$\kappa = \frac{2 + 16 C_H'}{C_m}, \qquad \delta = \frac{1}{8\kappa \|m\|} \wedge \frac{\log 2}{C_H'}, \qquad \gamma = \frac{\log 8}{\delta}. \tag{63}$$

The initial arrangement and the routing rules will depend on $\hat{X}^{0,n}$ and, in particular, on whether $e \cdot \hat{X}^{0,n} < -1$ or not. First consider the case where $e \cdot \hat{X}^{0,n} < -1$.

*Initial arrangement.* Recall that $X^n(0) = X^{0,n}$ is given and we have to specify $\Psi^n(0)$. We set

$$\Psi_{21}^n(0) = \lceil n^{5/8} \kappa \rceil, \tag{64}$$

$$\Psi_{11}^n(0) = \lceil (N_1^n - N_2^n + X_1^{0,n} + X_2^{0,n})/2 \rceil - \Psi_{21}^n(0), \tag{65}$$

$$\Psi_{12}^n(0) = \lfloor (N_2^n - N_1^n + X_1^{0,n} - X_2^{0,n})/2 \rfloor + \Psi_{21}^n(0), \tag{66}$$

$$\Psi_{22}^n(0) = X_2^{0,n} - \Psi_{21}^n(0). \tag{67}$$

Using (4) and (5), one verifies that

$$Y_1^n(0) = Y_2^n(0) = 0, \qquad |Z_1^n(0) - Z_2^n(0)| \leq 1. \tag{68}$$

*Routing.* The routing decisions can be based only on $n$, the value of $\Xi^n$ right before the arrival, and some auxiliary randomness that we have denoted by $R^n$. The way we use the randomness in the proposed policy is so as to split the customers of class 2 into two sub-classes, which we label as $\alpha$ and $\beta$. Upon the $k$th class-2 arrival, an independent biased coin is tossed according to which it is decided what sub-class the arrival belongs to. The bias of the coin is allowed to depend on $k$. We denote by $\alpha_n(k)$ the probability that the $k$th class-2 arrival is classified as a class-$\alpha$ customer. The value of $\alpha_n(k)$ is determined as

$$\alpha_n(k) = \left\{ \frac{\kappa(\gamma + \mu_{21})}{\lambda_2 n^{3/8}} \exp\left[\frac{\gamma(k-1)}{n\lambda_2}\right] \right\} \wedge 1. \tag{69}$$

Since we assume that $\hat{X}^{0,n} < -1$ and the difference between $Z_1^n(0)$ and $Z_2^n(0)$ is at most 1, it follows that the initial arrangement is such that there are free servers in both stations. Below, we shall describe the routing policy only as long as there are free servers in both stations. The description of the policy at other times is not important and will be left completely open (in fact, one of the main ingredients of the proof will be to show that the event that there are no free servers in one of the stations some time before $T$ has probability approaching zero as $n \to \infty$). The routing rules are as follows:

1. Class-1 customers are routed to the station with more free servers. More precisely, if a class-1 customer arrives at time $t$, it is instantaneously routed to station $j$, where

$$j = \begin{cases} 1, & \text{if } \hat{Z}_1^n(t-) > \hat{Z}_2^n(t-), \\ 2, & \text{otherwise.} \end{cases}$$

2. Class-$\alpha$ customers are routed to station 1 upon arrival.
3. Class-$\beta$ customers are routed to station 2 upon arrival.

Next, consider the case where the initial condition satisfies $e \cdot \hat{X}^{0,n} \geq -1$. Denote $r_n = e \cdot X^{0,n} - e \cdot N^n + \lceil n^{1/2} \rceil$. For the initial arrangement, we let $r_n$ class-1 customers be left in the queue, and let $\Psi_{ij}^n(0)$ be defined as in (64)–(67), except that we substitute $X^{0,n} - r_n$ for $X^{0,n}$. As a result, in place of the left part of (68) we will have $Y_1^n(0) = r_n$, $Y_2^n(0) = 0$. The routing is determined as follows. The $r_n$ class-1 customers that are initially put in the queue are kept in the queue. Rules 1–3 above apply for all the other customers. Upon the first arrival after the time $\tau_0^n$ when the number of free servers in the system first exceeds $r_n + \lceil n^{1/2} \rceil$, all the $r_n$ customers are moved into service: $\gamma_1^n$ into station 1 and $\gamma_2^n$ into station 2. Here, $\gamma_1^n = \lceil (Z_2^n(\tau_0^n-) - Z_1^n(\tau_0^n-) + r_n)/2 \rceil$, $\gamma_2^n = \lfloor (Z_1^n(\tau_0^n-) - Z_2^n(\tau_0^n-) + r_n)/2 \rfloor$. Clearly, $e \cdot \hat{X}^n(\tau_0^n) \leq -1$. Also, by the above choice for $\gamma_i^n$, one achieves that $Y_i^n(\tau_0^n) = 0$ and $|Z_1^n(\tau_0^n) - Z_2^n(\tau_0^n)| \leq 1$, a situation similar to (68). From the time $\tau_0^n$ on, the routing rules 1–3 above apply for all the new arrivals. Note that once again we have ignored the scenario that there are no free servers in one of the stations at some time before $T$, for reasons as mentioned before.

Here and for most of this section we shall assume that $e \cdot x < -1$. Since $\hat{X}^{0,n} \to x$, we have

$$e \cdot \hat{X}^{0,n} < -1 \tag{70}$$

for all large $n$ and, as a result, only the first part of the definition of the policy will be relevant until, at the end of the section, we return to the treatment of the case $e \cdot x \geq -1$.

Without loss of generality, assume that $T$ is a multiple of $\delta$, and set $\bar{k} = T/\delta$. Divide the time interval $[0,T]$ as follows: $I_0 = \{0\}$, $I_1 = (0,\delta], \ldots, I_{\bar{k}} = ((\bar{k}-1)\delta, \bar{k}\delta]$. Let $\bar{h}: [0,T] \to \mathbb{R}$ be defined as

$$\bar{h}(t) = 2^{-k} \qquad \text{for all } t \in I_k,\ k = 0, 1, \ldots, \bar{k}.$$

Define

$$\tau^n = \inf\left\{t \geq 0 : \hat{X}_e^n(t) \geq -\frac{\bar{h}(t)}{2}\right\} \wedge T. \tag{71}$$

Let $h' = \frac{1}{32} 2^{-\bar{k}} \equiv \frac{1}{32} \bar{h}(T)$. Let

$$\sigma^n = \inf\{t \geq 0 : \hat{Z}_1^n(t) \wedge \hat{Z}_2^n(t) \leq 4h'\} \tag{72}$$

and $\zeta^n = \tau^n \wedge \sigma^n$. By (70) and since $\hat{Z}_1^n(0)$ and $\hat{Z}_2^n(0)$ are nearly equal [cf. (68)], it is useful to note that

$$\hat{Z}_1^n(0) \wedge \hat{Z}_2^n(0) \geq 1/4. \tag{73}$$

Also, it is clear that there are idle servers in both stations up to time $\sigma^n$. Denote

$$D(a,p) = \{\xi \in \mathbb{R}^I : \|\xi\| \leq a,\ e \cdot \xi \leq -p\}, \qquad a > 0,\ p > 0.$$

Let $K_n = n^{1/8}$ and

$$D_k^n = D(e^{\gamma \delta k} K_n, 2^{-k}).$$

For $k = 0, 1, \ldots, \bar{k}$, denote

$$\Omega_k^n = \{\hat{X}^n(t) \in D_{k'}^n,\ t \in I_{k'},\ \text{for all } k' \leq k\} \cap \{\zeta^n > k\delta\}. \tag{74}$$

PROPOSITION 4.1. *Let the assumptions of Theorem 2.4 hold and assume also that $e \cdot x < -1$. Then for $k = 0, 1, \ldots, \bar{k}$, $P(\Omega_k^n) \to 1$.*

The above result implies that $P(\sigma^n \leq T) \to 0$. Since by construction of the policy, $Y^n(t) = 0$ for $t \leq \sigma^n$, this establishes Theorem 2.4 in this case.

In preparation for the proof of Proposition 4.1 we need some notation and preliminary results. Recall that we have denoted the only cycle by $c$, and that, according to our notation, $\Psi_c^n(t) \equiv \Psi_{21}^n(t)$ and $\hat{\Psi}_c^n(t) \equiv \hat{\Psi}_{21}^n(t)$. Let $A_\alpha^n(t)$ $[A_\beta^n(t)]$ denote the number of class-$\alpha$ (resp., class-$\beta$) customers that have arrived up to time $t$. Let

$$\hat{A}_\alpha^n(t) = n^{-1/2} A_\alpha^n(t) - K_n \int_0^t (\mu_{21} + \gamma) \psi_0(s) \, ds, \tag{75}$$

where

$$\psi_0(t) = \kappa e^{\gamma t}, \qquad t \geq 0. \tag{76}$$

A sequence of processes is said to be $C$-tight if it is tight and every subsequential limit has continuous sample paths with probability one. The following three lemmas will be proved later in this section.

LEMMA 4.1. *Under the assumptions of Proposition 4.1, one has*

$$|\hat{A}_\alpha^n|_T^* \to 0 \qquad \textit{in probability}, \tag{77}$$

$$\|\bar{X}^n - x^*\|_{T \wedge \sigma^n}^* \equiv n^{-1/2} \|\hat{X}^n\|_{T \wedge \sigma^n}^* \to 0 \qquad \textit{in probability}, \tag{78}$$

$$\|\bar{\Psi}^n - \psi^*\|_{T \wedge \sigma^n}^* \equiv n^{-1/2} \|\hat{\Psi}^n\|_{T \wedge \sigma^n}^* \to 0 \qquad \textit{in probability}. \tag{79}$$

*In addition, the processes $\hat{A}_i^n(\cdot \wedge \sigma^n)$, $\hat{S}_{ij}^n(\int_0^{\cdot \wedge \sigma^n} \bar{\Psi}_{ij}^n(s) \, ds)$ are $C$-tight. In particular, one has $|\hat{S}_{21}^n(\int_0^{\cdot} \bar{\Psi}_{21}^n(s) \, ds)|_{T \wedge \sigma^n}^* \to 0$ in probability.*

LEMMA 4.2. *Under the assumptions of Proposition 4.1, $|\hat{\Psi}_c^n - K_n \psi_0|_{T \wedge \sigma^n}^* \to 0$ in probability, where $\psi_0$ is as in (76).*

LEMMA 4.3. *Under the assumptions of Proposition 4.1, $\lim_{n \to \infty} P(\sigma^n \leq \tau^n) = 0$.*

The following will be used several times in the proofs of both Proposition 4.1 and the above lemmas. Recall that, for $t \leq \sigma^n$, $Y^n(t) = 0$. Hence, we have from (27)

$$\hat{X}^n(t) = \hat{X}^{0,n} + \hat{W}^n(t) + \int_0^t [H^n(s) + m^n \hat{\Psi}_c^n(s)] \, ds, \qquad t \leq \sigma^n, \tag{80}$$

where we have denoted

$$H^n(s) = H^n(\hat{X}^n(s), -\hat{Z}^n(s)). \tag{81}$$

Also, substituting zero for $\hat{Y}^n$ in (28), using (14) and positivity of $\hat{Z}_j^n$, we have, for all $n$ large,

$$\|\hat{Z}^n(t)\| \leq \|\hat{X}^n(t)\| + 1, \qquad t \leq \sigma^n. \tag{82}$$

Thus, by (25), (62), (81) and (82), for all $n$ large,

$$\|H^n(s)\| \vee \|e \cdot H^n(s)\| \leq C'_H(1 + \|\hat{X}^n(s)\|), \qquad t \leq \sigma^n. \tag{83}$$

PROOF OF PROPOSITION 4.1. We argue by induction on $k$. As the induction base, consider $k = 0$. By (70) and since $\hat{X}^{0,n}$ converge, we have $\hat{X}^n(0) \in D_0^n$ for all large $n$. By Lemma 4.3, it remains to show that $\tau^n > 0$ with probability approaching 1. This follows from (71) and the fact that the jumps of $\hat{X}^n$ are bounded by $n^{-1/2}$.

Next consider the induction step. Assuming that for a given $k < \bar{k}$ one has

$$\lim_{n \to \infty} P(\Omega_k^n) = 1, \tag{84}$$

we shall prove

$$\lim_{n \to \infty} P(\Omega_{k+1}^n) = 1. \tag{85}$$

In view of Lemma 4.3, it suffices to show

$$\lim_{n \to \infty} P(\tau^n \leq (k+1)\delta) = 0, \tag{86}$$

and

$$\lim_{n \to \infty} P(\|\hat{X}^n\|_{\delta(k+1) \wedge \sigma^n}^* \leq e^{\gamma\delta(k+1)} K_n) = 1. \tag{87}$$

To this end, note that from (80) we have

$$\hat{X}_e^n(t) = \hat{X}_e^{0,n} + \hat{W}_e^n(t) + \int_0^t [H_e^n(s) - C_m^n \hat{\Psi}_c^n(s)] ds, \qquad t \leq \sigma^n. \tag{88}$$

By (63), one verifies that the constants $\kappa$, $\delta$ and $\gamma$ satisfy

$$e^{C'_H \delta} \leq 2, \qquad 1 + C'_H e^{\gamma\delta} \leq \tfrac{1}{2} C_m \kappa, \tag{89}$$

$$4(1 + \|m\|\kappa\delta e^{\gamma\delta}) \leq e^{\gamma\delta}. \tag{90}$$

Denote

$$\widetilde{\Omega}_k^n = \Omega_k^n \cap \{\sigma^n > \tau^n\} \cap \{\tau^n \leq (k+1)\delta\}.$$

Let $h = \tfrac{1}{8} 2^{-k}$. On $\Omega_k^n$ (hence, on $\widetilde{\Omega}_k^n$) the following must hold:

$$\|\hat{X}^n(k\delta)\| \leq e^{\gamma\delta k} K_n, \qquad \hat{X}_e^n(k\delta) \leq -8h. \tag{91}$$

QUEUES IN HEAVY TRAFFIC 31

In addition, on $\widetilde{\Omega}_k^n$,

there exist $k\delta < s < t \leq (k+1)\delta \wedge \sigma^n$ such that
$$\hat{X}_e^n(s) \leq -7h, \qquad \hat{X}_e^n(t) \geq -5h.$$

Hence, in view of (88), the following must hold on $\widetilde{\Omega}_k^n$:

$$\begin{aligned}(92)\qquad 2h &\leq \hat{X}_e^n(t) - \hat{X}_e^n(s) \\ &= \hat{W}_e^n(t) - \hat{W}_e^n(s) + \int_s^t [H_e^n(u) - C_m^n \hat{\Psi}_c^n(u)]\, du.\end{aligned}$$

Denote $a_n = |\hat{\Psi}_c^n - K_n \psi_0|_{T \wedge \sigma^n}^*$, $b_n = \|\hat{W}^n\|_{T \wedge \sigma^n}^*$ and $d_{u,v}^n = \hat{W}_e^n(v) - \hat{W}_e^n(u)$. Let

$$(93)\qquad \hat{\Omega}_k^n = \Omega_k^n \cap \{C_m a_n \leq 2\} \cap \{2b_n + 2\delta \|m\| a_n + T \leq K_n\}.$$

By Lemmas 4.1–4.3 and (84),

$$(94)\qquad \lim_{n \to \infty} P(\hat{\Omega}_k^n) = 1.$$

On the event $\hat{\Omega}_k^n$, by (80) and (83), denoting
$$\Lambda^n = \|\hat{X}^n(k\delta)\| + 2b_n + \delta \|m^n\| \, |\hat{\Psi}_c^n|_{T \wedge \sigma^n}^* + T,$$

we have
$$\|\hat{X}^n(u)\| \leq \Lambda^n + \int_{k\delta}^u C_H' \|\hat{X}^n(u')\|\, du', \qquad u \leq (k+1)\delta \wedge \sigma^n.$$

By Gronwall's inequality and the first part of (89), this shows that $\|\hat{X}^n(u)\| \leq 2\Lambda^n$, for $u$ as above. Using this along with (93), and assuming in what follows that $n$ is so large that $\|m^n\| \leq 2\|m\|$ and $C_m^n \geq C_m/2$, we have on the event $\hat{\Omega}_k^n$

$$(95)\qquad \|\hat{X}^n(u)\| \leq 4K_n e^{\gamma \delta k}(1 + \|m\|\kappa \delta e^{\gamma \delta}) \leq K_n e^{\gamma \delta (k+1)},$$
$$u \leq (k+1)\delta \wedge \sigma^n,$$

where in the last inequality we used (90). Equations (94) and (95) establish (87). Next, combining (92) and (95), and using again (93), we have on the event $\widetilde{\Omega}_k^n \cap \hat{\Omega}_k^n$

$$2h \leq d_{s,t}^n + (C_H' + 1)(t - s) + K_n(t - s)e^{\gamma \delta k}[C_H' e^{\gamma \delta} - C_m \kappa/2].$$

Writing $C = C_H' + 1$ and using the second part of (89), we conclude that, for all large $n$, on the event $\widetilde{\Omega}_k^n \cap \hat{\Omega}_k^n$, there exist $s$ and $t$ such that $k\delta < s < t \leq (k+1)\delta \wedge \sigma^n$ and

$$d_{s,t}^n \geq 2h + (e^{\gamma \delta k} K_n - C)(t - s) \geq \begin{cases} \frac{1}{2} K_n^{1/2}, & t - s \geq K_n^{-1/2}, \\ 2h, & t - s < K_n^{-1/2}. \end{cases}$$

Hence, for all large $n$,

$$P(\widetilde{\Omega}_k^n) \leq P((\hat{\Omega}_k^n)^c) + P(\text{there exist } 0 \leq s < t \leq T \wedge \sigma^n \text{ such that } 2d_{s,t}^n \geq K_n^{1/2})$$
$$+ P(\text{there exist } 0 \leq s < t \leq T \wedge \sigma^n$$
$$\text{such that } t - s < K_n^{-1/2} \text{ and } d_{s,t}^n \geq 2h)$$
$$\leq P((\hat{\Omega}_k^n)^c) + P(4b_n \geq K_n^{1/2}) + P(w_{T \wedge \sigma^n}(\hat{W}_n, K_n^{-1/2}) \geq h).$$

By (94), the tightness of $b_n$ and the $C$-tightness statement in Lemma 4.1, we see that $P(\widetilde{\Omega}_k^n) \to 0$ as $n \to \infty$. By (84) and Lemma 4.3, this shows (86). □

Let a constant $r > 0$ be given and let $\theta^n = \inf\{t \in [0,\delta] : \hat{X}_e^n(t) \leq -r\}$. Let $\varepsilon$ be as in Theorem 2.4. The following consequence of the above proof will be useful.

COROLLARY 4.1. *Under the assumptions of Proposition 4.1,* $\lim_{n \to \infty} P(\theta^n > \varepsilon/2) = 0$. *In addition,* $\|\hat{X}^n\|_{\theta^n}^*$ *are tight.*

PROOF. With the notation from (74) and (93), define $\Omega^n = \{\theta^n > \varepsilon/2\} \cap \hat{\Omega}_1^n$. By Lemmas 4.1–4.3 and Proposition 4.1, $\lim_{n \to \infty} P(\hat{\Omega}_1^n) = 1$ and, therefore, it suffices to show that $P(\Omega^n) \to 0$. On $\Omega^n$ we have

$$-r \leq \hat{X}_e^n(\varepsilon) - \hat{X}_e^n(0) = \hat{W}_e^n(\varepsilon) - \hat{W}_e^n(0) + \int_0^\varepsilon [H_e^n(u) - C_m^n \hat{\Psi}_c^n(u)] \, du,$$

and an argument along the lines of the proof of Proposition 4.1 proves the first claim of the result.

Next we prove the tightness statement. Note that, by assumption, $\|\hat{X}^n(0)\| \leq C_0$ for some deterministic constant $C_0$ independent on $n$. Let $\varsigma^n = (2r + \|e\|C_0)e^{-\gamma\delta}K_n^{-1}$ and $\widetilde{\Omega^n} = \{\theta^n > \varsigma^n\} \cap \hat{\Omega}_1^n$. On $\widetilde{\Omega^n}$ we have (just as in the proof of Proposition 4.1)

$$-r \leq \hat{X}_n^e(\varsigma^n) = \hat{X}_e^n(0) + \hat{W}_e^n(\varsigma^n) + \int_0^{\varsigma^n} [H_e^n(u) - C_m^n \hat{\Psi}_c^n(u)] \, du$$
$$\leq \hat{X}_e^n(0) + \hat{W}_e^n(\varsigma^n) + (C_H' + 1)\varsigma^n - K_n \varsigma^n e^{\gamma\delta}$$
$$\leq \hat{W}_e^n(\varsigma^n) + (C_H' + 1)\varsigma^n - 2r.$$

As a result,

(96) $$\lim_{n \to \infty} P(\theta^n > \varsigma^n) = 0.$$

Let $\widetilde{C} = C_0 + 4r\|m\|\kappa + 2C_H r + 1$. On $\hat{\Omega}_1^n$, for $n$ so large that $\varsigma^n < \delta$, with the notation $a^n = |\Psi_c^n - K_n \psi_0|_{T \wedge \sigma^n}^*$, we have

$$\|\hat{X}^n\|_{\varsigma^n}^* \leq \|\hat{X}^n(0)\| + w_T(\hat{W}^n, \varsigma^n)$$
$$+ \|m^n\|\delta a^n + \|m^n\|\kappa K_n e^{\gamma \delta} \varsigma^n + C_H K_n e^{\gamma \delta} \varsigma^n$$
$$\leq w_T(\hat{W}^n, \varsigma^n) + 2\|m\|\delta a^n + \widetilde{C} - 1.$$

Since, by Lemmas 4.1 and 4.2, we have $\lim_{n \to \infty} P(w_T(\hat{W}^n, \varsigma^n) + 2\|m\|\delta a^n > 1) = 0$, we conclude that $\lim_{n \to \infty} P(\|\hat{X}^n\|_{\varsigma^n}^* \leq \widetilde{C}) = 1$. In view of (96), this shows that $\|\hat{X}^n\|_{\theta^n}^*$ are tight. $\square$

We turn to prove Lemmas 4.1–4.3.

PROOF OF LEMMA 4.1. We first prove (77). Let $R_k^n$ denote the indicator of the event that, in the $n$th system, the $k$th class-2 arrival was classified as an arrival of class $\alpha$. Then $ER_k^n = P(R_k^n = 1) = \alpha_n(k)$ given in (69). Since $A_\alpha^n(t)$ represents the number of class-2 arrivals up to time $t$ that were classified as class $\alpha$, we can write $A_\alpha^n(t) = \sum_{k=1}^{A^n(t)} R_k^n$. By (75), we have

$$\hat{A}_\alpha^n(t) = n^{-1/2} A_\alpha^n(t) - n^{1/8} \frac{\lambda_2 C}{\gamma} (e^{\gamma t} - 1),$$

where throughout this proof $C = (\mu_{21} + \gamma)\kappa/\lambda_2$. Fix $\varepsilon > 0$. Then $\hat{A}_\alpha^n(t) \geq \varepsilon$ if and only if $A_\alpha^n(t) \geq n^{1/2}\varepsilon + K_t^n$, where

$$K_t^n = n^{5/8} \frac{\lambda_2 C}{\gamma} (e^{\gamma t} - 1).$$

Therefore,

(97) $\qquad P[\text{there exists } t \leq T, \text{ such that } \hat{A}_\alpha^n(t) \geq \varepsilon] \leq p_1^n + p_2^n,$

where

$$p_1^n = P\left(\sup_{t \leq T} |A_2^n(t) - \lambda_2^n t| \geq n^{3/4}\right) \to 0 \qquad \text{as } n \to \infty$$

(98)
$$p_2^n = P\left(\sup_{t \leq T} |A_2^n(t) - \lambda_2^n t| < n^{3/4}, \text{ and there exists } t \leq T, \right.$$
$$\text{such that } \sum_{k=1}^{A^n(t)} R_k^n \geq n^{1/2}\varepsilon + K_t^n \right)$$
$$\leq P\left(\text{there exists } t \leq T, \text{ such that } \sum_{k=1}^{\beta(n,t)} R_k^n - K_t^n \geq n^{1/2}\varepsilon\right),$$

where $\beta(n,t) = \lfloor \lambda_2^n t + n^{3/4} \rfloor$. The convergence statement in (98) is due to the tightness of $\hat{A}_2^n$ [cf. (56)]. A direct calculation based on (69) shows

$$\sum_{k=1}^{\beta(n,t)} E(R_k^n) \leq K_t^n + C' n^{3/8},$$

for all large $n$, where $C'$ is a deterministic constant not depending on $n$. Hence,

$$p_2^n \leq P\left[\text{there exists } t \leq T \text{ such that } \sum_{k=1}^{\beta(n,t)} (R_k^n - ER_k^n) \geq n^{1/2}\varepsilon - C'n^{3/8}\right].$$

Denoting $m(l) = \sum_{k=1}^{l}(R_k^n - ER_k^n)$, we can write

$$p_2^n \leq P(|m|_{2\lambda^n T}^* \geq n^{1/2}\varepsilon/2) \leq 16 \frac{E[m(2\lambda_2^n T)^2]}{\varepsilon^2 n},$$

where in the last step we used Doob's inequality for the martingale $m$. In turn, $E[m(l)^2] = \sum_{k=1}^{l} \alpha_n(k)$, and substituting $2\lambda_2^n T$ for $l$, one finds that $p_2^n \to 0$. As a result, the l.h.s. of (97) converges to zero. A similar calculation, which we omit, shows that the probability that there exists $t \leq T$ such that $\hat{A}_\alpha^n(t) \leq -\varepsilon$ also converges to zero. Since $\varepsilon$ is arbitrary, (77) follows.

Exactly as in the proof of Theorem 2.3, we have

$$(99) \qquad n^{-1/2}\hat{W}^n \Rightarrow 0.$$

Since $\hat{X}^{0,n}$ converges, we have from (80), (81) and (83), for all large $n$,

$$(100) \quad \begin{aligned} \|\hat{X}^n(t)\| &\leq K_n + \|\hat{W}^n(t)\| \\ &\quad + C_H' \int_0^t (1 + \|\hat{X}^n(s)\|) \, ds + T\|m^n\| \|\hat{\Psi}_c^n\|_t^*, \qquad t \leq \sigma^n. \end{aligned}$$

Using the inequality

$$(101) \qquad \hat{\Psi}_c^n(t \wedge \sigma^n) \leq \hat{\Psi}_c^n(0) + n^{-1/2} A_\alpha^n(T)$$

and Gronwall's inequality, we have

$$(102) \quad \|\hat{X}^n(t)\| \leq C_1(K_n + \|\hat{W}^n(t)\| + n^{-1/2}A_\alpha^n(T) + 1), \qquad t \leq \sigma^n \wedge T,$$

for an appropriate constant $C_1$ not depending on $n$. It follows from (77) that

$$(103) \qquad \lim_{n \to \infty} P(A_\alpha^n(T) \geq n^{3/4}) = 0.$$

Combining (99), (102) and (103), we have (78). Using (4), (5), (11), (16)–(18) and the fact that $\hat{Y}^n(t) = 0$, we have, for all $t \leq \sigma^n$,

$$\hat{\Psi}_{11}^n(t) + \hat{\Psi}_{12}^n(t) = \hat{X}_1^n(t),$$
$$\hat{\Psi}_c^n(t) + \hat{\Psi}_{22}^n(t) = \hat{X}_2^n(t),$$
$$\hat{Z}_1^n(t) + \hat{\Psi}_{11}^n(t) + \hat{\Psi}_c^n(t) = \hat{N}_1^n.$$

These equations along with (82) show that

(104) $$\|\hat{\Psi}^n(t)\| \le C_2(1 + \|\hat{X}^n(t)\| + \hat{\Psi}_c^n(t)), \qquad t \le \sigma^n,$$

for a constant $C_2$ not depending on $n$. Combining (101), (103) and (78), we have (79). The result (79) implies that the processes $\bar{\Psi}_{ij}^n(\cdot \wedge \sigma^n)$ are tight and that every subsequential limit has continuous sample paths with probability one. Using this and the time change lemma [4] along the lines of the proof of Theorem 2.3 above proves $C$-tightness as claimed. Finally, consider the last claim in the statement of the lemma. Since $\psi_{21}^* = 0$, (79) shows that $|\bar{\Psi}^n|_{T \wedge \sigma^n}^* \to 0$ in probability. Using this and (56), the claim follows using again the time change lemma. $\square$

PROOF OF LEMMA 4.2. By construction of the policy, $B_{21}^n(t) = A_\alpha^n(t)$ for all $t \le \sigma^n$. Thus, by (7),

$$\Psi_{21}^n(t) = \Psi_{21}^n(0) + A_\alpha^n(t) - S_{21}^n\left(\int_0^t \Psi_{21}^n(s)\,ds\right).$$

Denote

(105) $$\varrho(t) = (\mu_{21} + \gamma)\psi_0(t), \qquad t \ge 0.$$

Recalling that $\hat{\Psi}_c^n(t) = n^{-1/2}\Psi_c^n(t) = n^{-1/2}\Psi_{21}^n(t)$, using (75), one checks by direct calculation that

$$\hat{\Psi}_c^n(t) = \hat{\Psi}_c^n(0) + n^{1/8}\int_0^t \varrho(s)\,ds - \mu_{21}^n \int_0^t \hat{\Psi}_c^n(s)\,ds + \widetilde{W}_0^n(t), \qquad t \le \sigma^n,$$

where

$$\widetilde{W}_0^n(t) = \hat{A}_\alpha^n(t) - \hat{S}_{21}^n\left(\int_0^t \bar{\Psi}_{21}^n(s)\,ds\right).$$

Note that $\psi_0(t) = \kappa + \int_0^t \varrho(s)\,ds - \mu_{21}\int_0^t \psi_0(s)\,ds$, and let $\psi_0^n$ be the unique solution to

(106) $$\psi_0^n(t) = \kappa + \int_0^t \varrho(s)\,ds - \mu_{21}^n \int_0^t \psi_0^n(s)\,ds.$$

Then, for $t \le \sigma^n$,

$$\hat{\Psi}_c^n(t) - n^{1/8}\psi_0^n(t) = \hat{\Psi}_c^n(0) - n^{1/8}\kappa$$
$$- \mu_{21}^n \int_0^t [\hat{\Psi}_c^n(s) - n^{1/8}\psi_0^n(s)]\,ds + \widetilde{W}_0^n(t).$$

By (64), $\hat{\Psi}_c^n(0) - n^{1/8}\kappa \to 0$. By Lemma 4.1, $|\widetilde{W}_0^n|_{T \wedge \sigma^n}^* \to 0$ in probability. An application of Gronwall's inequality therefore shows that

(107) $$|\hat{\Psi}_c^n - n^{1/8}\psi_0^n|_{T \wedge \sigma^n}^* \to 0 \qquad \text{in probability.}$$

Also, by (13) and (106), it is easy to see that $\delta^n(t) \leq C\varepsilon_n + C\int_0^t \delta^n(s)\,ds$ for $t \leq T$, where $\delta^n(t) = |\psi_0(t) - \psi_0^n(t)|$, $C$ is a constant and $n^{1/8}\varepsilon_n \to 0$. Hence, $n^{1/8}|\delta^n|_T^* \to 0$. Along with (107), this proves the lemma. □

PROOF OF LEMMA 4.3. By (73), the initial values for both $\hat{Z}_1^n$ and $\hat{Z}_2^n$ are greater than $8h'$. We thus have $\{\sigma^n \leq \tau^n\} \subseteq \Omega_1^n \cup \Omega_2^n$, where, for $i = 1, 2$,

$$\Omega_i^n = \Big\{\text{there exist } u_1 < u_2 \leq \zeta^n \text{ such that}$$

$$\sup_{u \in [u_1, u_2]} \hat{Z}_i^n(u) \leq 7h', \hat{Z}_i^n(u_1) \geq 6h', \hat{Z}_i^n(u_2) \leq 5h'\Big\}.$$

*Step* 1. We show that $P(\Omega_1^n) \to 0$. Let $u_1$ and $u_2$ be as specified in the expression for $\Omega_1^n$. Note that the size of the jumps of $\hat{X}^n$ is bounded by $n^{-1/2}$. Hence, for $u \leq \tau^n$, we have $\hat{X}_e^n(u) \leq -32h' + n^{-1/2}$. This, combined with (14) and (28), and the fact that $Y_i^n(u) = 0$ for $u \leq \sigma^n$ imply that $\hat{Z}_1^n(u) < \hat{Z}_2^n(u)$ for $u \in [u_1, u_2]$. According to our definition of the policy, all class-1 customers are routed to station 2 during the period $[u_1, u_2]$ and only class-$\alpha$ customers are routed to station 1 during this period. As a result, $B_{11}^n(u_2) - B_{11}^n(u_1) = 0$ and $B_{21}^n(u_2) - B_{21}^n(u_1) = A_\alpha^n(u_2) - A_\alpha^n(u_1)$. By (7), we have

$$\sum_{i=1,2}(\Psi_{i1}^n(u_2) - \Psi_{i1}^n(u_1)) = A_\alpha^n(u_2) - A_\alpha^n(u_1) - D_1^n(u_1, u_2),$$

where

$$D_j^n(t_1, t_2) := \sum_{i=1,2}\left[S_{ij}^n\left(\int_0^{t_2}\Psi_{ij}^n(s)\,ds\right) - S_{ij}^n\left(\int_0^{t_1}\Psi_{ij}^n(s)\,ds\right)\right], \qquad j = 1, 2.$$

Note that $D_j^n(t_1, t_2)$ represents the number of departures from station $j$ during $(t_1, t_2]$. Using (5) and (14), one has on $\Omega_1^n$ that

(108) $$h'n^{1/2} \leq A_\alpha^n(u_2) - A_\alpha^n(u_1) - D_1^n(u_1, u_2).$$

Denoting, for $j = 1, 2$,

$$\widetilde{W}_j^n(t) := \sum_{i=1,2}\hat{S}_{ij}^n\left(\int_0^t \bar{\Psi}_{ij}^n(s)\,ds\right),$$

(109) $$\widetilde{S}_j^n(t_1, t_2) = \widetilde{W}_j^n(t_2) - \widetilde{W}_j^n(t_1)$$

$$+ \sum_{i=1,2}\left[\mu_{ij}^n\int_{t_1}^{t_2}\hat{\Psi}_{ij}^n(s)\,ds + n^{1/2}\mu_{ij}^n\psi_{ij}^*(t_2 - t_1)\right]$$

and

$$\widetilde{A}_1^n(t_1, t_2) = \hat{A}_\alpha^n(t_2) - \hat{A}_\alpha^n(t_1) + n^{5/8} \int_{t_1}^{t_2} \varrho(s)\, ds,$$

where $\varrho$ is as in (105), one checks by direct calculation that (108) can be written as $h' \leq \widetilde{A}_1^n(u_1, u_2) - \widetilde{S}_1^n(u_1, u_2)$. Note that $\varrho$ is bounded above. Moreover, $\mu_{11}^n \psi_{11}^*$ is bounded below by a positive constant. Using (104) and denoting $\Delta = u_2 - u_1$ we therefore have, for all large $n$,

(110) $$h' \leq 2|\hat{A}_\alpha^n|_T^* + w_{T \wedge \sigma^n}(\widetilde{W}_1^n, \Delta) - n^{1/2} \Lambda_1^n \Delta,$$

where

$$\Lambda_1^n = C_1 - C_2 n^{-1/2}(1 + \|\hat{X}^n\|_{T \wedge \sigma^n}^* + |\hat{\Psi}_c^n|_{T \wedge \sigma^n}^*),$$

and $C_1 > 0$ and $C_2$ are suitable constants. We conclude that $P(\Omega_1^n) \leq P(\Omega_{1,1}^n) + P(\Omega_{1,2}^n)$, where

$$\Omega_{1,1}^n = \{\text{there exists } \Delta \in (0, n^{-1/4}] \text{ such that (110) holds}\},$$

$$\Omega_{1,2}^n = \{\text{there exists } \Delta \in (n^{-1/4}, T] \text{ such that (110) holds}\}.$$

By Lemmas 4.1 and 4.2, $\Lambda_1^n \to C_1$ and $|\hat{A}_\alpha^n|_T^* \to 0$ in probability, and $\widetilde{W}_1^n(\cdot \wedge \sigma^n)$ are $C$-tight. This shows that $P(\Omega_{1,1}^n) \to 0$. On $\Omega_{1,2}^n$, if $\Lambda_1^n \geq 0$, then $n^{1/4} \Lambda_1^n \leq \varepsilon_n$ must hold, where $\varepsilon_n = 2|\hat{A}_\alpha^n|_T^* + 2|\widetilde{W}_1^n|_{T \wedge \sigma^n}^*$. Hence, $P(\Omega_{1,2}^n) \leq P(\Lambda_1^n < 0) + P(\Lambda_1^n \leq n^{-1/4} \varepsilon_n) \to 0$. This shows that $P(\Omega_1^n) \to 0$.

*Step* 2. We next show that $P(\Omega_2^n) \to 0$. Letting $u_1$ and $u_2$ be as in the expression for $\Omega_2^n$ and arguing as before, one obtains that $\hat{Z}_1^n(u) > \hat{Z}_2^n(u)$ for $u \in [u_1, u_2]$ and, consequently, that all class-1 customers are routed to station 1 during this period. Analogously to (108), we find that on $\Omega_2^n$

$$h' n^{1/2} \leq A_\beta^n(u_2) - A_\beta^n(u_1) - D_2^n(u_1, u_2).$$

Using the inequality $A_\beta^n(u_2) - A_\beta^n(u_1) \leq A_2^n(u_2) - A_2^n(u_1)$ and some direct calculation, we deduce from the above that

(111) $$h' \leq \widetilde{A}_2^n(u_1, u_2) - \widetilde{S}_2^n(u_1, u_2),$$

where $\widetilde{S}_2^n$ is as in (109) and

$$\widetilde{A}_2^n(t_1, t_2) = \hat{A}_2^n(t_2) - \hat{A}_2^n(t_1) + \lambda_2^n n^{-1/2}(t_2 - t_1).$$

Note that the r.h.s. of (111) contains the term $\gamma_n n^{1/2}(u_2 - u_1)$, where $\gamma_n = n^{-1}\lambda_2^n - \sum_{i=1,2} \mu_{i2}^n \psi_{i2}^*$. We claim that $\gamma_n$ is bounded above by a negative constant for all $n$ large. Indeed, recall that $\xi_{21}^* = 0$ and note that from (10) we have $\lambda_2 = \bar{\mu}_{22}\xi_{22}^*$. With (9) and (11), this shows $\lambda_2 = \mu_{22}\psi_{22}^* < \sum_{i=1,2} \mu_{i2}\psi_{i2}^*$. The claim regarding $\gamma_n$ thus follows from (8).

Along the lines of step 1, we obtain, instead of (110),

$$h' \leq w_{T \wedge \sigma^n}(\hat{A}_2^n - \widetilde{W}_2^n, \Delta) - n^{1/2}\Lambda_2^n \Delta,$$

where $\Delta = u_2 - u_1$,

$$\Lambda_2^n = C_3 - C_4 n^{-1/2}(1 + \|\hat{X}^n\|_{T \wedge \sigma^n}^* + |\hat{\Psi}_c^n|_{T \wedge \sigma^n}^*),$$

and $C_3 > 0$ and $C_4$ are constants. The rest of the argument for showing $P(\Omega_2^n) \to 0$ is similar to that in step 1 and is omitted. We conclude that $P(\sigma^n \leq \tau^n) \to 0$. $\square$

PROOF OF THEOREM 2.4. In view of Proposition 4.1, it only remains to treat the case where $e \cdot x \geq -1$. We can split the sequence of systems into two subsequences according as $e \cdot \hat{X}^{0,n} < -1$ or not, and since the result is already proved for the subsequence on which $e \cdot \hat{X}^{0,n} < -1$, we will assume, without loss of generality, that $e \cdot \hat{X}^{0,n} \geq -1$ for all $n$. As a result, the second part of the definition of the proposed policy applies. In the first part of the proof we used the fact that the system starts with no queues. In the current situation we argue that, even if the system starts with a queue, it is brought quickly to zero. Recall the notation $\tau_0^n$ from definition of the policy. Note that, prior to time $\tau_0^n$, there are $r_n$ class-1 customers that are kept in the queue and that, as far as all other customers are concerned, the system behaves exactly as in the first part of the definition. As a result, we can use Corollary 4.1 with $r = 1 + \sup_n n^{-1/2} r_n < \infty$, and it follows that $P(\tau_0^n \leq \varepsilon/2)$ converges to 1. At time $\tau_0^n$ the system is in a state very similar to that in which a system satisfying $e \cdot \hat{X}^{0,n} < -1$ is at time zero, in the sense that $e \cdot \hat{X}^n(\tau_0^n) < -1$ and $|Z_1^n(\tau_0^n) - Z_2^n(\tau_0^n)| \leq 1$. A review of the proof of Proposition 4.1, replacing the initial values of all processes by their values at time $\tau_0^n$, shows that, with probability approaching 1, $Y^n(t)$ is kept zero for $t \in [\varepsilon, T]$. As in Section 3, the only remaining issue is that $\hat{X}^n(\tau_0^n)$ must be shown to be tight. This is indeed the case by Corollary 4.1. $\square$

## APPENDIX A

LEMMA A.1. (i) *All nodes $\mathcal{I} \cup \mathcal{J}$ of $\mathcal{G}_b$ are connected through the edges of $\mathcal{G}_b$. (ii) Every nonbasic activity $(i,j) \in \mathcal{E}_{nb}$ belongs to a simple cycle.*

PROOF. For (i), consider a node of $\mathcal{G}$, $i \in \mathcal{I}$. By (10), $\sum_j \bar{\mu}_{ij} \xi_{ij}^* > 0$ and, therefore, there exists $j \in \mathcal{J}$ such that $\xi_{ij}^* > 0$. This shows that $(i,j) \in \mathcal{E}_b$ and thus, $i$ must be a node of $\mathcal{G}_b$. A similar argument holds for a node $j \in \mathcal{J}$ observing that $\sum_i \xi_{ij}^* = 1$ by the heavy traffic condition. Item (ii) follows since $\mathcal{G}_b$ is a tree. $\square$

PROOF OF THEOREM 2.1. We first show that the relations (3)–(6), when rescaled appropriately, imply the following:

$$(112) \quad \hat{X}_i^n(t) = \hat{X}_i^{0,n} + \hat{W}_i^n(t) - \sum_{j \in \mathcal{J}} \mu_{ij}^n \int_0^t \hat{\Psi}_{ij}^n(s)\,ds, \qquad i \in \mathcal{I}, j \in \mathcal{J},$$

$$(113) \qquad \hat{Y}_i^n(t) + \sum_{j \in \mathcal{J}} \hat{\Psi}_{ij}^n(t) = \hat{X}_i^n(t), \qquad i \in \mathcal{I},$$

$$(114) \qquad \hat{Z}_j^n(t) + \sum_{i \in \mathcal{I}} \hat{\Psi}_{ij}^n(t) = \hat{N}_j^n, \qquad j \in \mathcal{J},$$

$$(115) \quad \hat{Y}_i(t) \geq 0, \qquad \hat{Z}_j(t) \geq 0, \qquad \hat{\Psi}_c^n(t) \geq 0, \qquad i \in \mathcal{I}, j \in \mathcal{J}, c \in \mathcal{C}, t \geq 0.$$

To this end, note that by (15), (16) and (19),

$$\begin{aligned}
\hat{X}_i^n(t) &= n^{-1/2}(X_i^n(t) - nx_i^*) \\
&= \hat{X}_i^{0,n} + \hat{A}_i^n(t) - \sum_{j \in \mathcal{J}} \hat{S}_{ij}^n\left(\int_0^t \bar{\Psi}_{ij}^n(s)\,ds\right) \\
&\quad + n^{-1/2}\left[\lambda_i^n t - \sum_{j \in \mathcal{J}} n\mu_{ij}^n \int_0^t \bar{\Psi}_{ij}^n(s)\,ds\right].
\end{aligned}$$

By (13), (17) and (18), the last term in the above display is equal to

$$-\sum_{j \in \mathcal{J}} \mu_{ij}^n \int_0^t \hat{\Psi}_{ij}^n(s)\,ds + \hat{\lambda}_i^n t + n^{1/2}\left[\lambda_i - \sum_{j \in \mathcal{J}} \mu_{ij}^n \psi_{ij}^*\right]t.$$

Using (10), (11) and (13),

$$n^{1/2}\left[\lambda_i - \sum_{j \in \mathcal{J}} \mu_{ij}^n \psi_{ij}^*\right]t = -\sum_{j \in \mathcal{J}} \hat{\mu}_{ij}^n \psi_{ij}^* t.$$

Combined with (26), this establishes (112) above. Equations (113) and (114) and the first two inequalities in (115) directly follow from (4), (5), (6), (11), (14), (16), (17) and (18). For the last inequality in (115), recall that, for $c \in \mathcal{C}$, $\hat{\Psi}_c^n = \hat{\Psi}_{ij}^n$, where $(i,j) = \sigma^{-1}(c)$, and that, for all $(i,j)$ of this form, $\psi_{ij}^* = 0$.

Equations (28)–(29) follow from (113)–(115). As for (27), it suffices to prove the identity

$$(116) \qquad \begin{aligned} \hat{\Psi}_{ij}^n(t) &= G_{ij}(\hat{X}^n(t) - \hat{Y}^n(t), \hat{N}^n \\ &\quad - \hat{Z}^n(t)) - \sum_{c \in \mathcal{C}:(i,j) \in c} s(c,i,j)\hat{\Psi}_c^n(t), \end{aligned}$$

for all $(i,j) \in \mathcal{E}$, $t \geq 0$. Indeed, with (21), the above implies the identity

$$-\sum_{j \in \mathcal{J}} \mu_{ij}^n \hat{\Psi}_{ij}^n(t) = -\sum_{j \in \mathcal{J}} \mu_{ij}^n G_{ij}(\hat{X}^n(t) - \hat{Y}^n(t), \hat{N}^n - \hat{Z}^n(t))$$
$$+ \sum_{c \in \mathcal{C}} m_{i,c}^n \hat{\Psi}_c^n(t),$$

which, along with (25) and (112), establishes (27). In order to show (116), define

(117) $$\Upsilon_{ij}^n(t) = \hat{\Psi}_{ij}^n(t) + \sum_{c \in \mathcal{C} \,:\, (i,j) \in c} s(c,i,j) \hat{\Psi}_c^n(t), \qquad (i,j) \in \mathcal{E}.$$

As follows from the uniqueness statement succeeding (22), in order to prove (116), it is enough to show that $\{\Upsilon_{ij}^n\}$ satisfy the system of equations (22) with $a = \hat{X}^n(t) - \hat{Y}^n(t)$ and $b = \hat{N}^n - \hat{Z}^n(t)$. To this end, note that if $(i,j)$ is a nonbasic activity, then there is exactly one simple cycle $c$ for which $(i,j)$ is an edge and one has $s(c,i,j) = -1$ and, by definition of $\hat{\Psi}_c^n$, $\Psi_c^n = \Psi_{ij}^n$. This shows that $\Upsilon_{ij}^n = 0$ for $(i,j) \in \mathcal{E}_{nb}$. Next, note that, for a given $c \in \mathcal{C}$ and a node $i \in \mathcal{I}$ of $c$, there are exactly two edges of $c$ of the form $(i,j)$. By the construction of directions, for these two values of $j$, the signifiers $s(c,i,j)$ must have opposite signs. As a result, $\sum_{c \in \mathcal{C}} \sum_{j \,:\, (i,j) \in c} s(c,i,j) \hat{\Psi}_c^n(t) = 0$. Changing the order of summation and using an analogous argument for summation over $i$, we obtain

(118)
$$\sum_{j \in \mathcal{J}} \sum_{c \in \mathcal{C} \,:\, (i,j) \in c} s(c,i,j) \hat{\Psi}_c^n(t) = 0, \qquad i \in \mathcal{I},$$
$$\sum_{i \in \mathcal{I}} \sum_{c \in \mathcal{C} \,:\, (i,j) \in c} s(c,i,j) \hat{\Psi}_c^n(t) = 0, \qquad j \in \mathcal{J}.$$

Equations (113), (114) and (118) now show that $\{\Upsilon_{ij}^n\}$ satisfy (22) with the values $a$ and $b$ mentioned above. This proves (116) and, in turn, the relation (27), and completes the proof of the theorem. $\square$

PROOF OF LEMMA 3.1. Let $X^n(t)$, $Y^n(t)$, $Z^n(t)$ and $\{\Psi_c^n(t)\}$ be given, satisfying all requirements in the statement of the lemma. Define $\{\Psi_{ij}^n(t)\}$ according to (42). Equations (2), (4) and (5) follow directly from the definition of $G$ and $\{s(c,i,j)\}$. As for (6), it remains only to prove that $\Psi_{ij}^n(t) \geq 0$ for every $(i,j) \in \mathcal{E}_b$. Along the lines of Lemma 3 of [2], note that $\psi^* = G(x^*, \nu)$, as follows from (11). Using linearity of the map $G$ on the domain $D_G$ and the bound assumed on $\Psi_c^n$,

$$\Psi_{ij}^n(t) \geq G_{ij}(nx^*, n\nu)$$
$$+ G_{ij}(X^n(t) - nx^* - Y^n(t), N^n - n\nu - Z^n(t)) - C_1 a_0 n,$$

where $C_1$ is a constant independent of $n$, $t$ and the choice of $a_0$ and, therefore,

$$\Psi_{ij}^n(t) \geq (\psi_{ij}^* - C_1 a_0) n \\ - C_2(\|X^n(t) - nx^* - Y^n(t)\| \vee \|N^n - n\nu - Z^n(t)\|),$$

where $C_2$ is a constant depending only on the map $G$. With the assumed bound on $\|X^n - nx^*\| \vee \|Y^n\| \vee \|Z^n\|$, using also (14), we obtain $\Psi_{ij}^n(t) \geq (\psi_{ij}^* - C_3 a_0) n - C_4 n^{1/2}$ for appropriate constants $C_3$, $C_4$ independent of $n$, $t$ and $a_0$. Since $\psi_{ij}^* > 0$ for all $(i,j) \in \mathcal{E}_b$, it follows that $a_0 > 0$ may be chosen so that, for all large $n$, $\Psi_{ij}^n(t) \geq 0$. □

SKETCH OF PROOF OF PROPOSITION 2.1. We simplify the notation by writing $A^n$, $\bar{A}^n$, $X^n$ and, respectively, $\bar{X}^n$ for $A_{OL}^n$, $\bar{A}_{OL}^n$, $X_{OL}^n$, $\bar{X}_{OL}^n$, and so on (there will be no confusion). Note that an analogous result to (56) holds and, in particular, for every $t$,

$$(119) \qquad n^{-1} A_i^n(t) \to t \lambda_i' \qquad \text{in probability},$$

as $n \to \infty$. Fix $t \geq 1$. By (112), we have

$$\bar{X}_i^n(t) = \bar{X}_i^n(0) + \bar{A}_i^n(t) - \sum_j \mu_{ij} \int_0^t \bar{\Psi}_{ij}^n(s) \, ds + \Delta_i^n(t),$$

where

$$(120) \quad \Delta_i^n(t) = \sum_j \left[ n^{-1/2} \hat{S}_{ij}^n \left( \int_0^t \bar{\Psi}_{ij}^n(s) \, ds \right) + (\mu_{ij}^n - \mu_{ij}) \int_0^t \bar{\Psi}_{ij}^n(s) \, ds \right].$$

By (5) and (6), $\bar{\Psi}_{ij}^n(t)$ is bounded. An argument as in step 2 in the proof of Theorem 2.3 shows that the first term in (120) converges to zero in probability. Since $\mu_{ij}^n \to \mu_{ij}$, so does the second term and, in turn, $\Delta_i^n(t) \to 0$ in probability. It can be shown that there exists a constant $\delta > 0$ depending only on $(\lambda, \lambda', \bar{\mu})$, such that, for any sub-stochastic matrix $(\xi_{ij})$,

$$(121) \qquad \max_i \left[ \lambda_i' - \sum_j \bar{\mu}_{ij} \xi_{ij} \right] \geq \delta.$$

Let $\Omega^n = \{\text{for all } i \in \mathcal{I}, \ |\Delta_i^n(t)| + |t^{-1} \bar{A}_i^n(t) - \lambda_i'| \leq \delta/3\}$. By (119), we have $P(\Omega^n) \to 1$. Define

$$\xi_{ij}^n = \frac{1}{t} \int_0^t \frac{\bar{\Psi}_{ij}^n(s)}{\nu_j} \, ds, \qquad \bar{\xi}_{ij}^n = \frac{\xi_{ij}^n}{\sum_{i'} \xi_{i'j}^n}.$$

It can be shown by (5) and (8) that $\lambda_i' - \sum_j \bar{\mu}_{ij} \xi_{ij}^n \geq \lambda_i' - \sum_j \bar{\mu}_{ij} \bar{\xi}_{ij}^n - \delta/3$, for all $n$ large and all $i$. Since $\bar{\xi}$ is sub-stochastic, (121) implies that, for some $i$, we must have $\lambda_i' - \sum_j \bar{\mu}_{ij} \bar{\xi}_{ij}^n \geq \delta$, hence, $\lambda_i' - \sum_j \bar{\mu}_{ij} \xi_{ij}^n \geq 2\delta/3$, for some $i$.

This can be used to show $\max_i \bar{X}_i^n(t) \geq \delta t/3$ on $\Omega^n$. Due to the relation $\bar{Y}_i^n(t) + \sum_j \bar{\Psi}_{ij}^n(t) = \bar{X}_i^n(t)$, we have $\bar{Y}_i^n(t) \geq \bar{X}_i^n(t) - C$, for all $i$, where $C$ does not depend on $n$ or $t$. Hence, $\max_i \bar{Y}_i^n(t) \geq \delta t/3 - C$. We conclude that on $\Omega^n$, for all $n$ large, $\sum_i Y_i^n(t) \geq [\delta t/3 - C]n$, and the result follows. $\square$

R. ATAR
DEPARTMENT OF ELECTRICAL ENGINEERING
TECHNION
HAIFA 32000
ISRAEL
E-MAIL: atar@ee.technion.ac.il

A. MANDELBAUM
G. SHAIKHET
DEPARTMENT OF INDUSTRIAL
 ENGINEERING AND MANAGEMENT
TECHNION
HAIFA 32000
ISRAEL